\documentclass[preprint,onefignum,onetabnum]{siamart190516}

\usepackage{lipsum}
\usepackage{amsfonts}
\usepackage{amscd}
\usepackage{mathabx}
\usepackage{graphicx}
\usepackage{epstopdf}
\usepackage{algorithmic}
\ifpdf
  \DeclareGraphicsExtensions{.eps,.pdf,.png,.jpg}
\else
  \DeclareGraphicsExtensions{.eps}
\fi

\usepackage{tikz}
\usepackage{tikz-cd}
\usetikzlibrary{arrows,cd,external,decorations.pathreplacing,calligraphy}
\usepackage{tikz-3dplot}
\usepackage{pgfplots}
\pgfplotsset{compat=1.15}
\usepgfplotslibrary{groupplots}
\usetikzlibrary{pgfplots.groupplots}
\usetikzlibrary{plotmarks}
\usetikzlibrary{calc}
\usepgfplotslibrary{external}
\usepackage{siunitx}
\usepackage{color}
\usepackage{enumerate}
\usepackage{subfigure}
\usepackage[titletoc,title]{appendix}

\usepackage{cleveref}
\usepackage{cancel}

\graphicspath{{./}{figures/}{Figures/}{../plots_numerical_results}}

\newsiamremark{remark}{Remark}
\newsiamremark{hypothesis}{Hypothesis}
\crefname{hypothesis}{Hypothesis}{Hypotheses}
\newsiamthm{claim}{Claim}

\headers{High Order geometric methods with splines}{B. Kapidani and R. V\'azquez}

\title{
  High order geometric methods with splines: \\an analysis of discrete Hodge--star operators\thanks{Submitted to the editors \today \funding{
  This work was supported by the Swiss National Science Foundation via the project HOGAEMS n.200021\_188589.}}
  }

\author{Bernard Kapidani\thanks{Institute of Mathematics, Ecole Polytechnique F\'ed\'erale de Lausanne, Station 8, 1015 Lausanne, Switzerland
  (\email{bernard.kapidani@epfl.ch}, \email{rafael.vazquez@epfl.ch}).}
\and Rafael V\'azquez\footnotemark[2] \thanks{Istituto di Matematica Applicata e Tecnologie Informatiche ``E.\ Magenes'' del CNR, via Ferrata 5, 27100 Pavia, Italy}.}

\usepackage{amsopn}

\newcommand{\Rd}{\color{red}}
\newcommand{\Bd}{\color{blue}}

\newcommand{\bmixed}{\beta}

\newcommand{\Xhath}[1]{\,\hspace{-1pt}\widehat{X}^{#1}_h\hspace{-1pt}\,}
\newcommand{\Xhathb}[1]{\,\hspace{-1pt}\widehat{X}^{#1}_{h,0}\hspace{-1pt}\,}

\newcommand{\Xhat}{\hat X}

\newcommand{\XhathT}[1]{\,\hspace{-1pt}\widehat{\tilde{X}}^{#1}_h\hspace{-1pt}\,}
\newcommand{\XhT}[1]{\,\hspace{-1pt}{\widetilde{X}}^{#1}_h\hspace{-1pt}\,}

\let\hat\widehat
\let\tilde\widetilde

\def\SC#1{\ifnum `#1<`1 {\Gd \hat {X}^0_h} \else \ifnum  `#1<`2
  \Bd{\hat X^1_h} \else \ifnum  `#1<`3 \Rd{\hat X^2_h} \else \Md{\hat X^3_h} \fi  \fi  \fi}
\def\XC#1{\ifnum `#1<`1 {\Gd {X}^0_h} \else \ifnum  `#1<`2
  \Bd{ X^1_h} \else \ifnum  `#1<`3 \Rd{ X^2_h} \else \Md{ X^3_h} \fi  \fi  \fi}
\def\XCt#1{\ifnum `#1<`1 {\Gd \tilde {X}^0_h} \else \ifnum  `#1<`2
  \Bd{ \tilde X^1_h} \else \ifnum  `#1<`3 \Rd{ \tilde X^2_h} \else \Md{ \tilde X^3_h} \fi  \fi  \fi}

\ifpdf
\hypersetup{
  pdftitle={High order geometric methods with splines: an analysis of discrete Hodge--star operators},
  pdfauthor={Bernard Kapidani and Rafael V\'azquez}
}
\fi

\begin{document}

\maketitle





\begin{abstract}
A new kind of spline geometric method approach is presented. Its main ingredient is the use of well established spline spaces forming a discrete de Rham complex to construct a primal sequence $\{X^k_h\}^n_{k=0}$, starting from splines of degree $p$, and a dual sequence $\{\tilde{X}^k_h\}_{k=0}^n$, starting from splines of degree $p-1$. By imposing homogeneous boundary conditions to the spaces of the primal sequence, the two sequences can be isomorphically mapped into one another.
Within this setup, many familiar second order partial differential equations can be finally accommodated by explicitly constructing appropriate discrete versions of constitutive relations, called Hodge--star operators.
Several alternatives based on both global and local projection operators between spline spaces will be proposed.
The appeal of the approach with respect to similar published methods is twofold: firstly, it exhibits high order convergence. Secondly, it does not rely on the geometric realization of any (topologically) dual mesh. Several numerical examples in various space dimensions will be employed to validate the central ideas of the proposed approach and compare its features with the standard Galerkin approach in Isogeometric Analysis.
\end{abstract}

\begin{keywords}
  Spline complex, Hodge--star, de Rham complex, Geometric method, isogeometric analysis
\end{keywords}

\begin{AMS}
  65D07, 65N12, 65N30 
\end{AMS}

\section{Introduction} \label{sec:introduction}
The expressive power of the language of exterior calculus and differential forms has long been recognized among physicists~\cite{flanders_differential_1989}. 
Its use allows for the neat distinction between topological and metric properties. 
Conservation laws in dynamical systems are for example succintly written using exterior derivatives, wedge products and Hodge--star operators. The former two are used to express system invariants (e.g. energy), elegantly written as duality pairings of differential forms whose order adds up to the dimension of the ambient space. The latter object encodes metric information of the domain, expressed either under the form of phenomenological constitutive equations or of the metric tensor when the physical problem is set on a manifold domain. 
In the field of numerical analysis, these tools have been particularly successful in the case of electromagnetics, where the symmetry and invariance properties of Maxwell's equations are directly embedded in the classical de Rham complex of differential forms. In this respect, the work of Bossavit~\cite{BOS88,BOS98b} and Tonti~\cite{tonti_finite_2001} has been pivotal in making the computational science community aware of the possible benefits of using differential forms, while the seminal works by Hiptmair~\cite{HIP99c,HIP02a} gave a solid foundation for their numerical analysis. 

One of the main issues in developing methods based on discrete differential forms is the definition of a discrete version of the Hodge--star operators \cite{HIP99c}, and depending on how these discrete operators are defined, we can separate these methods into two broad families:
\begin{enumerate}[i)]
\item Galerkin--Hodge approaches, in which the discrete Hodge--star weakly satisfies the same equation as its continuous counterpart, when tested against locally supported polynomials approximations of the involved differential forms.
\item Geometric Hodge approaches, in which 
the discrete Hodge--star operator is defined from the primal to the dual sequence (and vice versa) by exploiting the geometric relationships between two meshes of the same physical domain.
\end{enumerate}

Among the many works that fall into the realm of Galerkin methods we would like to emphasize the ones of Bossavit and Kettunen~\cite{TKB99,Bossavit_Kettunen_2000}, who gave the first interpretation of the discrete Hodge--star operator in finite elements, the review paper by Hiptmair~\cite{HIP02a}, and more recently the introduction of the fruitful field of finite element exterior calculus by Arnold, Falk and Winther~\cite{AFW06,AFW-2}.

On the more applied side, there has been has been plenty of effort dedicated to the development of numerical methods in the second family. 
Methods based on discrete exterior calculus~\cite{DEC-arxiv,Hirani-darcy,hirani}, which explicitly construct a second dual mesh by taking centroids of elements of the first mesh and connecting them with edges~\cite{hirani,lipnikov_mimetic_2014}, are very popular.
Methods based on the integral formulation of Maxwell's equations such as the finite integration technique (FIT)~\cite{ClThWe99,Clemens-Weiland}, or the cell method~\cite{tonti_finite_2001,m_marrone_computational_2001,codecasa_explicit_2008,codecasa_novel_2018} also fall into this second family. 
In general, for the differential form approximations, 
one proceeds by attaching discrete degrees of freedom to geometric entities of the mesh (vertices, edges, faces and volumes), and by repeatedly using an adapted version of the generalized Stokes theorem, which amounts to building incidence matrices for the geometric entities. If circumcentric dual meshes are used, discrete Hodge--star operators become diagonal matrices relating the geometric entities from the primal and the dual mesh, leading to very fast methods, even though in general the error decreases at most linearly upon mesh refinement.

The two families are not entirely disjoint, since the explicit construction of a dual mesh has been also considered in the finite element setting. In particular, to build a dual de Rham complex with stable pairing in two dimensions, one needs to work on a dual grid obtained by barycentric refinement~\cite{BuCh07}, mainly used as a preconditioner to speed up boundary element formulations~\cite{Andriulli_2008} or as a Lagrange multiplier for mesh coupling~\cite{niemimaki_structure-preserving_2016}.
In addition to being currently limited to two-dimensional manifolds, the construction does not lend itself to higher orders of approximation.
One of the present authors has recently introduced an extension of the cell method for the two dimensional Maxwell initial value problem, achieving arbitrary order of convergence and block diagonal mass matrices~\cite{kapidani_arbitrary-order_2021}, but since the approach is based on non-conforming spaces~\cite{kapidani_time-domain_2020}, sequence properties are not trivially preserved, and differential operators are not clearly related to the topology of an underlying mesh.
In any case, the explicit geometric construction of a topological dual mesh can be a cumbersome procedure which adds one layer of complexity to the algorithm and usually introduces more unknowns. 

Isogeometric discrete differential forms were recently introduced to extend finite element exterior calculus to the emerging framework of Isogeometric Analysis (IGA). A discrete de Rham complex was first constructed and analyzed for tensor-product B-splines \cite{Buffa_Sangalli_Vazquez,BRSV11}, though generalizations based on analysis-suitable T-splines \cite{BSV14}, locally refined B-splines \cite{Johannessen15} and hierarchical B-splines \cite{evans_hierarchical_2020} have also been proposed aiming at local refinement. These spaces have been applied in the Galerkin framework for the discretization of Maxwell's equations \cite{ratnani,Corno20161} also in the context of plasma physics \cite{kraus_kormann_morrison_sonnendrucker_2017}, and for the development of pointwise divergence free methods for incompressible fluid flow \cite{Buffa_deFalco_Sangalli,EvHu12,EvHu12-2,EvHu12-3,Van17}. While most of the previous works were based on vector fields, the first attempts to fully exploit the framework of differential forms for B-splines were made in \cite{ratnani,Back_2014} and \cite{Hiemstra20141444}, with discrete Hodge--star operators fitting the Galerkin--Hodge approach.
Recently, Gerritsma and co-authors have introduced a new approach in \cite{jain_construction_2020,Gerritsma_curlcurl_2020}, in which the primal space is used as its dual, and the discretization of the Hodge--star operator is replaced by the discretization of the codifferential, which is defined as the composition of the exterior derivative with two Hodge--star operators, on its left and its right.

In this work we introduce a new method based on isogeometric differential forms which falls into the second family presented above, in the sense that it is based on the explicit construction of a dual spline complex. In contrast to the above mentioned methods of the same family, a dual grid is not explicitly built, and the dual complex is simply defined by a change in polynomial degree, with the same construction used in \cite{Buffa_2020aa,kapidani2021treecotree} for stable mortar coupling between non-conforming meshes. Thanks to the high continuity of splines, the exterior derivative is rigorously defined in both sequences, and given by incidence matrices of a Cartesian grid \cite{ratnani,BSV14}, and the dimension of pairing spaces from the two sequences is always equal. The complete definition of the method requires, apart from the dual complex, a set of discrete Hodge--star operators, and we will analyze the properties of three different ones: the first two are a direct adaption of \cite{HIP99c} to our setting, while the third one is a local projector belonging to the family of quasi-interpolants defined by Lee, Lyche and M{\o}rken in \cite{LLM01}. We will prove that, when applied to elliptic problems, the method attains high order of convergence with any of the three operators, although one order of convergence is lost with respect to a Galerkin method based on the primal complex. We also present numerical evidence showing that when applied to Maxwell eigenvalue problem, the method is devoid of spurious solutions.

An outline of the paper is as follows.  In Section~\ref{sec:derham} we introduce the mathematical notation and define the de Rham complex of differential forms. In Section~\ref{sec:psplinecomplex} we present the B-spline complex of isogeometric discrete differential forms and discuss some of its properties. We introduce in Section~\ref{sec:dsplinecomplex} the dual spline complex by applying a change of the degree, and define the three different discrete Hodge--star operators that we will study. In Section~\ref{sec:erroranalysis} we provide some error estimates for the application of the resulting numerical approach to an elliptic model problem. In Section~\ref{sec:tests} we numerically validate convergence and high order approximation properties. Finally, in Section~\ref{Sec:conclusion} we summarise the main features of the present contribution and discuss future efforts in improving the new method.

\section{Differential forms} \label{sec:derham}

In the following we will give a brief recap of concepts related to differential forms and exterior calculus. The covered material is by no means original and thorough treatment of the topic can be found in \cite{AFW-2,HIP99c} and references therein.

\subsection{Alternating forms}
Let us denote with $\text{Alt}^k \mathbb{R}^n$, for $0 \le k \le n$, the space of alternating $k$--forms, linear maps $\omega$ : $\bigtimes_{i=1}^k \mathbb{R}^n \mapsto \mathbb{R}$ that assign to each $k$--tuple of vectors a real number, and that are linear in each argument and reverse sign for any odd permutation of the arguments. The dimension of the space of alternating $k$--forms in $n$--dimensional Euclidean space is given by
\begin{equation*}
    \text{dim}(\text{Alt}^k\mathbb{R}^n) = {n \choose k} ={n \choose n-k}, 
\end{equation*}
\noindent from which we can consistently set $\text{Alt}^0 \mathbb{R}^n := \mathbb{R}$.
Accordingly, we recognize that $\text{dim}\left(\text{Alt}^1 \mathbb{R}^n\right) = \text{dim}\left(\mathbb{R}^n\right)$ and introduce 1--forms 
through the linear map $dx_{\alpha}: \mathbb{R}^n \rightarrow \mathbb{R}$ 
which maps a vector to its $\alpha^{\textup{th}}$ coordinate. 
One constructs higher order alternating forms using 1--forms as the basic building block. To do so, the exterior (or wedge) product must be introduced. This is an associative operator $\wedge:\,\text{Alt}^k \mathbb{R}^n \times \text{Alt}^l \mathbb{R}^n \mapsto \text{Alt}^{k+l} \mathbb{R}^n$. For two 1--forms, it holds
\begin{equation}dx_{\alpha_i} \wedge dx_{\alpha_j} = -dx_{\alpha_j} \wedge dx_{\alpha_i},\label{eq:skew-sim}\end{equation}
where $1\leq i,j \leq n$, with the obvious implication $dx_{\alpha_i} \wedge dx_{\alpha_i} = 0$. In general the wedge product between alternating forms of higher degree is built upon the one between 1--forms, and takes the form
\begin{equation*}
w^k = \sum_{\boldsymbol \alpha \in \mathcal{I}_k} c_{\boldsymbol \alpha} \; dx_{\alpha_1} \wedge dx_{\alpha_2} \wedge \cdots \wedge dx_{\alpha_k},
\end{equation*}
where $w^k\in\text{Alt}^k \mathbb{R}^n$, $c_{\boldsymbol \alpha} \in \mathbb{R}$ and $\boldsymbol \alpha = (\alpha_1, \alpha_2, \ldots, \alpha_k) \in \mathcal{I}_k$, the set of ordered multi-indices
\begin{equation*}
\mathcal{I}_k = \left\{ \boldsymbol \alpha = (\alpha_1, \alpha_2, \ldots, \alpha_k): 1 \leq \alpha_1 < \alpha_2 < \ldots < \alpha_k \leq n \right\}.
\end{equation*}

In plain words, the object $dx_{\alpha_1} \wedge dx_{\alpha_2} \wedge \dots \wedge dx_{\alpha_k}$, often called a decomposable $k$--form, is an element in the basis of the space of alternating $k$--forms. Here it is assumed that $\alpha_1 < \alpha_2 < \dots < \alpha_k$, but the assumption is by no means necessary, and in fact it suffices to require
\begin{equation*}
dx_{\alpha_1} \wedge dx_{\alpha_2} \wedge \cdots \wedge dx_{\alpha_k} = \textup{sgn}(\sigma) dx_{\sigma(\alpha_1)} \wedge dx_{\sigma(\alpha_2)} \wedge \ldots \wedge dx_{\sigma(\alpha_k)},
\end{equation*}
\noindent where $\sigma: \left\{\alpha_1, \alpha_2, \ldots, \alpha_k\right\} \rightarrow \left\{\alpha_1, \alpha_2, \ldots, \alpha_k\right\}$ is an arbitrary permutation of the set $\left\{\alpha_1, \alpha_2, \ldots, \alpha_k\right\}$ and $\textup{sgn}(\sigma)$ is its sign (+1 if there is an even number of pairs $\alpha_i < \alpha_j$ such that $\sigma(\alpha_i) > \sigma(\alpha_j)$ and -1 otherwise).  This generalizes (\ref{eq:skew-sim}) to higher dimensions and implies
\begin{equation*}
dx_{\alpha_1} \wedge dx_{\alpha_2} \wedge \cdots \wedge dx_{\alpha_k} = 0, \; \text{ if $\alpha_i = \alpha_j$ for some $i \neq j$,}
\end{equation*}
\noindent as required by the alternating property. 
The space $\text{Alt}^k\mathbb{R}^n$ for $k>n$ is empty, while the space of alternating $n$--forms is entirely generated by the volume form 
\begin{equation*}
dV := dx_1 \wedge dx_2 \wedge \cdots \wedge dx_n.
\end{equation*}

\subsection{Differential forms, de Rham diagrams and the Hodge--star operator} \label{sec:derham-hodge}
Alternating forms are the building block of more general objects, namely the linear combination of alternating forms with more general functions (supported on a bounded subset of $\mathbb{R}^n$) as coefficients, called differential forms.
Let $D \subset \mathbb{R}^n$ a given open domain, that for simplicity we will assume to be contractible, and define a differential $0$--form as a function $\omega^0$: $D \mapsto \mathbb{R} =\text{Alt}^0\mathbb{R}^n$. 
Generally, for $k \ge 1$, a differential $k$--form will be a function $\omega^k$: $D\mapsto\text{Alt}^k\mathbb{R}^n$, which can be written as the linear combination, with differential $0$--forms as coefficients, of decomposable alternating $k$--forms:
\begin{equation*}
\omega^k = \sum_{\boldsymbol \alpha \in \mathcal{I}_k} \omega_{\boldsymbol \alpha} dx_{\alpha_1} \wedge dx_{\alpha_2} \wedge \cdots \wedge dx_{\alpha_k},
\end{equation*}
where $\omega_{\boldsymbol \alpha}: D \rightarrow \mathbb{R}$ is a differential 0--form,  and again $\boldsymbol \alpha = (\alpha_1, \alpha_2, \ldots, \alpha_k)$ is a multi-index. We denote the space of smooth differential $k$--forms with $\Lambda^k(D)$.

Accordingly, the notion of wedge product between a differential $k$--form and a differential $l$--form easily follows from the one involving alternating forms:
\begin{equation*}
\omega^k \wedge \eta^l := \sum_{\boldsymbol \alpha \in \mathcal{I}_k} \sum_{\boldsymbol \beta \in \mathcal{I}_l} \omega_{\boldsymbol \alpha} \eta_{\boldsymbol \beta}\, dx_{\alpha_1} \wedge \ldots \wedge dx_{\alpha_k} \wedge dx_{\beta_1} \wedge \ldots \wedge dx_{\beta_l},
\end{equation*}
and the result is a differential $(k+l)$--form. 
The introduction of functions as coefficients opens up the possibility of performing so-called exterior calculus on forms in addition to exterior algebra. 
For example, assuming some $0$--form $\omega^0$ to be smooth enough allows us to introduce its exterior derivative, defined as
\begin{equation*}
d^0\omega^0 := d\omega^0  = \sum_{\alpha=1}^n \frac{\partial \omega^0}{\partial x_{\alpha}} dx_\alpha. 
\end{equation*}
The above definition sheds light on a further connection between differential $0$--forms and $1$--forms since ${\partial \omega^0}/{\partial x_{\alpha}}$ is also a differential $0$--form, therefore implying $d\omega^0 \in \Lambda^1(D)$.
In general, making use of the exterior algebra machinery previously introduced, the exterior derivative operator for a differential $k$--form is then defined as the $(k+1)$--form
\begin{equation*}
d^k \omega^k = \sum_{\boldsymbol \alpha \in \mathcal{I}_k} d^0\omega_{\boldsymbol \alpha}^0 \wedge dx_{\alpha_1} \wedge dx_{\alpha_2} \wedge \cdots \wedge dx_{\alpha_k},
\end{equation*}
where $d^0\omega_{\boldsymbol \alpha}^0$ is the exterior derivative of the differential $0$--form $\omega_{\boldsymbol \alpha}^0$. From here onwards we will simply denote the exterior derivative by $d$ and when there is no confusion on the order of the differential form. 
An important property of the exterior derivative is 
 that $d \circ d \omega = 0$ for any differential $k$--form $\omega$.
 
In general, considering smooth functions is too restrictive, since we are interested in weak solutions of partial differential equations. Nevertheless, starting from the space of smooth differential $k$--forms $\Lambda^k(D)$, and similarly to what happens within more classic functional analytic approaches, we shall define the weighted $L^2$-inner product
\begin{equation*}
    \left( \omega^k, \eta^k \right)_{L_\gamma^2 \Lambda^k(D)} := 
    \sum_{\boldsymbol \alpha \in \mathcal{I}_k} \int_D \omega_{\boldsymbol\alpha}^0 \,\gamma\, \eta_{\boldsymbol\alpha}^0
    dx_{\alpha_1} \wedge dx_{\alpha_2} \wedge \cdots \wedge dx_{\alpha_n}, 
\end{equation*}
\noindent where we have added as a weight the parameter $\gamma\in L^\infty (D)$. This is in general a Riemannian metric tensor, but for the sake of clarity in the presentation we settle for a positive and uniformly bounded scalar valued parameter. When $\gamma$ is equal to one, the inner product reduces to the standard $L^2$-inner product, and we denote it by $(\cdot, \cdot)_{L^2 \Lambda^k(D)}$. 
By taking the completion of $\Lambda^k(D)$ with respect to this inner product we construct the Hilbert space denoted by $L^2\Lambda^k(D)$.

Now, let us define
\begin{equation*}
H\Lambda^k(D) = \left\{ \omega \in L^2\Lambda^k(D) : d^k \omega \in L^2\Lambda^{k+1}(D) \right\},
\end{equation*}
which is also a Hilbert space equipped with the inner-product
\begin{equation*}
\left( \omega, \eta \right)_{H\Lambda^k(D)} = \left( \omega, \eta \right)_{L^2 \Lambda^k(D)} + \left( d^k \omega, d^k \eta \right)_{L^2 \Lambda^{k+1}(D)}.
\end{equation*}
Similarly, we denote by $H^s \Lambda^k(\Omega)$ the space of differential forms such that all its coefficients $\omega_{\boldsymbol{\alpha}}$ belong to $H^s(\Omega)$ \cite{AFW06}.

With the above notation established, recalling the vanishing property of any composition between two exterior derivatives, a de Rham complex of differential forms can be constructed:
\begin{equation} \label{eq:derhamcomplex}
\begin{CD}
H\Lambda^0(D) @>d^0>> H\Lambda^1(D) @>d^1>> H\Lambda^2(D) @>d^2>> \ldots @>d^{n-1}>> H\Lambda^n(D),
\end{CD}
\end{equation}
\noindent where, for the last space in the sequence to be well defined, it is understood that the exterior derivative of an $n$--form identically vanishes. We refer to the above as the $L^2$ de Rham complex of differential forms. 
Furthermore, we will also need to define de Rham complexes of differential forms subject to homogeneous or vanishing boundary conditions. Such complexes are built from spaces of differential $k$--forms with compact support in exactly the same manner as in \eqref{eq:derhamcomplex}. We distinguish them by using the zero subscript and they build the sequence
\begin{equation} \label{eq:derhamcomplex_homogeneous}
\begin{CD}
H_0\Lambda^0(D) @>d^0>> H_0\Lambda^1(D) @>d^1>> H_0\Lambda^2(D) @>d^2>> \ldots @>d^{n-1}>> H_0\Lambda^n(D),
\end{CD}
\end{equation}
\noindent 
\noindent for which a fundamental result 
is given by the integration by parts formula 
\begin{equation} \label{eq:int_by_parts}
    \int_D d\omega \wedge \eta = 
(-1)^k \int_D \omega \wedge d \eta, \quad \text{ for } \omega \in H_0 \Lambda^{k-1}(D), \eta \in H\Lambda^{n-k} (D), 
\end{equation}
\noindent for $0 < k \leq n$. 
In the more general case of non-vanishing boundary conditions, an integral computed on $\partial D$ must be added to the right hand-side.
We remark that the assumption on the domain $D$ being contractible makes the de Rham complexes introduced in (\ref{eq:derhamcomplex}) and (\ref{eq:derhamcomplex_homogeneous}) exact, i.e. $Im(d^k) = Ker(d^{k+1})$, where we omit for the sake of brevity the obvious spaces and operators which start and end the sequences from and onto the trivial space.

The Hodge--star operator $\star_{\gamma}^k:L^2 \Lambda^k(D) \rightarrow L^2 \Lambda^{n-k}(D)$ is a linear operator from $k$--forms into twisted $(n-k)$--forms, which is defined by 
\begin{equation} \label{eq:def_hodge_diff}
    \left( \eta^k, \omega^k \right)_{L_\gamma^2 \Lambda^k(D)} = \int_D \eta^k \wedge \star_{\gamma}^k(\omega^k),
    \quad \forall \eta^k \in L^2 \Lambda^k(D), 
\end{equation}
\noindent where we note that the Hodge--star operator acquires a subscript related to the material parameter. 
For continuous $k$--forms the Hodge--star operator is invertible, with $\star^k_\gamma \circ \star^{n-k}_{1/\gamma} = (-1)^{k(n-k)}$~\cite{HIP99c}. This important property implies that 
\begin{equation} \label{eq:def_hodge_diff_inv}
(\eta^{n-k}, \star_\gamma^k (\omega^{k}))_{L^2_{1/\gamma}\Lambda^{n-k}(\Omega)} = (-1)^{k(n-k)} \int_\Omega \eta^{n-k} \wedge \omega^{k}, \quad \forall \eta^{n-k} \in L^2 \Lambda^{n-k}(\Omega).
\end{equation}
\noindent 
The fundamental equivalence of the Hodge--star operator with an $L^2$-inner product will be a guiding principle for the explicit construction of their discrete approximation.

\subsection{Parametric domain, physical domain, and pullback operators}
\label{subsec:pullback}
As customary in IGA, we will make extensive use of the open set $\hat \Omega = (0,1)^n$, which we refer to as the parametric domain, where we will first define discrete spaces of differential forms and only afterwards these will be mapped to a physical domain $\Omega \subset \mathbb{R}^n$ using a parametric mapping $\textbf{F}: \hat \Omega \rightarrow \Omega$. For the sake of readability we will then define the notation
\begin{equation*}
\hat{X}^k := H\Lambda^k(\hat\Omega), \quad X^k := H\Lambda^k(\Omega),
\end{equation*}
where we again add the $0$ subscript to denote subspaces with vanishing boundary conditions. We thus relate differential $k$--forms in the parametric domain to differential $k$--forms in the physical domain using a set of pullback operators $\iota^k: X^k \rightarrow \hat{X}^k$. The pullback operator for differential $0$--forms takes the form $\iota^0(\omega) = \omega \circ \textbf{F}$ for all $\omega \in X^0$, while the pullback operator for differential $n$--forms takes the form $\iota^n(\omega) = \det(D\textbf{F}) \left(\omega \circ \textbf{F}\right)$ for all $\omega \in X^n$ where $D\textbf{F}$ is the Jacobian matrix of the mapping $\textbf{F}$.  
An important property of pullbacks is that they commute with the exterior derivative.

\subsection{Model problem}
We will focus our analysis on the elliptic problem of finding a $k$--form $\omega \in H_0 \Lambda^k(\Omega)$ such that, for a given right-hand side $f \in L^2 \Lambda^k$, 
within the notation introduced above
\begin{equation} \label{eq:continuousproblem}
(-1)^{k+1} d^{n-k-1} \star_\gamma^{k+1} (d^k \omega) + \star^k_\beta \omega = \star^k f,
\end{equation}
\noindent holds. Although further regularity is needed for the correct definition of the second exterior derivative, for simplicity we keep the same regularity required for the variational formulation of the problem analyzed in Section~\ref{sec:erroranalysis}.
In the numerical results, we will also study the eigenvalue problem of finding $\omega \in H_0 \Lambda^k(\Omega)$ and $\tau \in \mathbb{R}$ such that 
\begin{equation} \label{eq:evp_with_forms}
(-1)^{k+1} d^{n-k-1} \star_\gamma^{k+1} (d^k \omega) = \tau \star^k_\beta \omega.
\end{equation}
We are particularly interested in Maxwell eigenvalue problem, which means that $k=1$, and the material properties are the magnetic reluctivity $\gamma \equiv \nu$, and the electric permittivity $\beta \equiv \varepsilon$.

\section{The B-spline complex of discrete differential forms} \label{sec:psplinecomplex}

In this section, we recall the tensor-product B-spline complex of isogeometric discrete differential forms, which we refer to henceforth as simply the B-spline complex \cite{BRSV11}. This complex is built using the tensor-product nature of multivariate B-splines as its name implies.  To fix notation, we first present univariate B-splines and multivariate B-splines, highlighting some of their properties which come in handy for the present work, before introducing the B-spline complex in both the parametric and physical domain. The notation is based mainly on \cite{evans_hierarchical_2020}.

\subsection{Univariate B-splines} \label{sec:univariate}
Let $p$ denote the polynomial degree of the univariate B-splines, and let $m$ denote the dimension of the space. To define the B-spline basis functions, we first introduce a $p$-open knot vector $\Xi = \{\xi_{1}, \ldots , \xi_{m+p+1} \}$, where
\begin{equation*}
0 = \xi_{1} = \ldots = \xi_{p+1} < \xi_{p+2} \le \ldots \le \xi_m < \xi_{m+1} = \ldots = \xi_{m+p+1} = 1.
\end{equation*}
Using the well known Cox-de Boor formula, we can then define B-spline basis functions $\left\{B_{i,p}\right\}_{i=1}^m$.  We denote by $S_p(\Xi)$ the space they span, which is the space of piecewise polynomials of degree $p$, with the number of continuous derivatives at each knot $\xi_i$ given by $p-r_i$, where $r_i$ is the multiplicity of the knot.
We can actually define the $i^{\textup{th}}$ B-spline basis function $B_{i,p}$ as the unique B-spline basis function defined from the Cox-de Boor algorithm and the local knot vector $\{\xi_i, \ldots, \xi_{i+p+1}\}$. Its support is given by the interval $[\xi_i,\xi_{i+p+1}]$.

Assuming the multiplicity of each of the internal knots is less or equal to $p$ (i.e., the B-spline functions are at least continuous), the derivative of a B-spline belonging to the space  $S_{p}(\Xi)$ belongs to the space $S_{p-1}(\Xi')$, where the new knot vector $\Xi' = \{\xi_2, \ldots, \xi_{m+p}\}$ is defined from $\Xi$ by removing the first and last repeated knots. Relatedly, we will also make use of the normalized basis functions $D_{i,p-1} = \frac{p}{\xi_{i+p} - \xi_i}B_{i,p-1}$, known as the Curry-Schoenberg basis: it holds, for the derivative of a B-spline, that 
\begin{equation*} 
B'_{i,p}(\zeta) = D_{i,p-1}(\zeta) - D_{i+1,p-1}(\zeta)
\end{equation*}
for $i = 1, \ldots, m$, where we use the convention $D_{1,p-1}(\zeta) = D_{m+1,p-1}(\zeta) = 0$ for any $\zeta \in (0,1)$. 

\subsection{Multivariate tensor-product B-splines}
Multivariate B-splines are constructed from a tensor-product of univariate B-splines.  Namely, given a vector of polynomial degrees ${\bf p} = (p_1, \ldots, p_n)$ and a set of open knot vectors ${\boldsymbol \Xi} = \left\{ \Xi_\ell \right\}_{\ell=1}^n$, we define the multivariate B-spline basis as
\[
{\cal B}_{\bf p}(\boldsymbol \Xi) := \{ B_{\bf i,p}({\boldsymbol \zeta}) = B_{i_1,p_1}(\zeta_1) \ldots B_{i_n,p_n}(\zeta_n) \} \; \textup{ for ${\boldsymbol \zeta} \in (0,1)^n$},
\]
and the corresponding space they span as
\begin{equation*}
S_{\bf p}(\boldsymbol \Xi) \equiv S_{p_1,\ldots,p_n}(\Xi_1,\ldots,\Xi_n) = \otimes^n_{\ell=1} S_{p_\ell}\left(\Xi_\ell\right) = {\rm span} ({\cal B}_{\bf p}(\boldsymbol \Xi)),
\end{equation*}
where we have introduced the shorthand $S_{\bf p}(\boldsymbol \Xi)$ to ease notation. Considering the knot vectors without repetitions we obtain a Cartesian mesh of the domain $(0,1)^n$. We will denote the mesh size by $h$, and assume shape regularity of the mesh, which implies local quasi-uniformity.

For the definition of the B-spline complex it will be needed to use tensor-products of spaces with mixed degree, combining the standard univariate spaces with the spaces of derivatives, in particular using the Curry-Schoenberg spline basis as done in \cite{ratnani}. We follow the notation in \cite{evans_hierarchical_2020} and introduce, for a given multi-index $\boldsymbol \alpha = \left(\alpha_1,\alpha_2,\ldots,\alpha_k\right)$, the space
\begin{equation*}
S_{\bf p}(\boldsymbol \Xi;\boldsymbol \alpha) = \otimes^n_{\ell=1} S_{p_\ell},
\end{equation*}
with the definition
\begin{equation*}
S_{p_\ell} \equiv S_{p_\ell}(\boldsymbol \Xi; \boldsymbol \alpha) = 
\left\{
\begin{array}{cl}
S_{p_\ell-1}(\Xi'_\ell) & \textup{ if } \ell = \alpha_j \textup{ for some } j \in \left\{1, \ldots, k \right\}, \\
S_{p_\ell}(\Xi_\ell) & \textup{ otherwise.}
\end{array} \right.
\end{equation*}
A basis of this space will be chosen by suitable combinations of standard B-splines for degree $p_\ell$, and Curry-Schoenberg splines for degree $p_\ell - 1$. More precisely, we define

\[
{\cal B}_{\bf p}(\boldsymbol \Xi; \boldsymbol \alpha) := \{\bmixed_{{\bf i}, {\bf p}}({\boldsymbol \zeta}) = \bmixed_{i_1,p_1}(\zeta_1) \ldots \bmixed_{i_n,p_n}(\zeta_n) \} \; \textup{ for ${\boldsymbol \zeta} \in (0,1)^n$},
\]
where
\[
\bmixed_{i_\ell,p_\ell} \equiv \bmixed_{i_\ell,p_\ell}(\boldsymbol \alpha) := 
\left \{
\begin{array}{cl}
D_{i_\ell,p_\ell-1} & \textup{ if } \ell = \alpha_j \textup{ for some } j \in \left\{1, \ldots, k \right\}, \\
B_{i_\ell,p_\ell} & \textup{ otherwise, }
\end{array}
\right.
\]
and it is immediate to see that $S_{\bf p}(\boldsymbol \Xi;\boldsymbol \alpha) = \mathrm{span} ({\cal B}_{\bf p}(\boldsymbol \Xi; \boldsymbol \alpha))$.

\subsection{The B-spline complex}\label{sec:spline-complex}
We are now able to define the B-spline complex of isogeometric discrete differential forms.  We first define it in the parametric domain before mapping it onto the physical domain.  To begin, let us assume that we are given a set of polynomial degrees $\left\{p_\ell\right\}_{\ell=1}^n$ and open knot vectors $\left\{ \Xi_\ell \right\}_{\ell=1}^n$, and let us further assume that multiplicity of the internal knots of $\Xi_\ell$ is never greater than $p_\ell$ for $\ell = 1, \ldots, n$, which implies that the functions in $S_{\bf p}(\boldsymbol \Xi)$ are continuous.  Then, we define the tensor-product spline space of isogeometric discrete differential $0$--forms as $\Xhath{0} := S_{\bf p}(\boldsymbol \Xi)$, and denote the corresponding basis by ${\cal B}^0 := {\cal B}_{\bf p}(\boldsymbol \Xi)$. For $k>0$, we define the basis of B-splines differential $k$--forms as 
\begin{equation*}
{\cal B}^k = \bigcup_{\boldsymbol \alpha \in {\cal I}_k} \left\{ \bmixed_{{\bf i}, {\bf p}} dx_{\alpha_1} \wedge \ldots \wedge dx_{\alpha_k} : \, \bmixed_{{\bf i}, {\bf p}} \in {\cal B}_{\bf p}(\boldsymbol \Xi;\boldsymbol \alpha) \right\},
\end{equation*}
and we denote the space they span, which is the tensor-product spline space of isogeometric discrete differential $k$--forms, by
\begin{equation*}
\Xhath{k} = \left\{ \hat{\omega}^h = \sum_{\boldsymbol \alpha \in \mathcal{I}_k} \hat{\omega}^h_{\boldsymbol \alpha} dx_{\alpha_1} \wedge \ldots \wedge dx_{\alpha_k} : \, \hat{\omega}^h_{\boldsymbol \alpha} \in S_{\bf p}(\boldsymbol \Xi;\boldsymbol \alpha), \hspace{5pt} \forall \boldsymbol \alpha \in \mathcal{I}_k \right\} = {\rm span}({\cal B}^k).
\end{equation*}
Then, it is easily shown that the above spaces constitute a discrete de Rham complex of the form
\begin{equation*} 
\begin{CD}
\Xhath{0} @>d^0>> \Xhath{1} @>d^1>> \ldots @>d^{n-1}>> \Xhath{n},
\end{CD}
\end{equation*}
which we call the B-spline complex.
To deal with spaces subject to vanishing boundary conditions, we define $\Xhathb{k} = \Xhath{k} \cap \Xhat^k_0$. It can again be easily shown that these spaces constitute a discrete de Rham complex of the form
\begin{equation} \label{eq:cd-iga-tp-bc}
\begin{CD}
\Xhathb{0} @>d^0>> \Xhathb{1} @>d^1>> \ldots @>d^{n-1}>> \Xhathb{n},
\end{CD}
\end{equation}
The discrete spaces in the physical domain are defined using the pullback operators defined in Section~\ref{subsec:pullback}, more precisely
\begin{equation*}
X_h^k := \left\{ \omega : \iota^k(\omega) \in \Xhath{k} \right\}, \quad
X_{h,0}^k := \left\{ \omega : \iota^k(\omega) \in \Xhathb{k} \right\}
\end{equation*}
for all integers $0 \leq k \leq n$. Using the definition of the pullback operators yields discrete de Rham complexes in the physical domain. In order to ensure good approximation properties we make the usual assumption from IGA that the mapping ${\bf F}$ is defined from the first space of the diagram $\hat{X}^0_h$, eventually using rational splines.

The meaning of the B-spline complex and the notation are more evident with the example of the three-dimensional setting, when using the isomorphisms between differential forms and scalar or vector proxy fields. In this case we find
\begin{align}
\Xhath{0} = & S_{\bf p}(\boldsymbol \Xi) \nonumber = S_{p_1,p_2,p_3}(\Xi_1,\Xi_2,\Xi_3), \nonumber \\
\Xhath{1} \simeq & S_{\bf p}\left(\boldsymbol \Xi; 1\right) \times S_{\bf p}\left(\boldsymbol \Xi; 2\right) \times S_{\bf p}\left(\boldsymbol \Xi; 3\right) \nonumber \\
= & S_{p_1-1,p_2,p_3}(\Xi'_1,\Xi_2,\Xi_3) \times S_{p_1,p_2-1,p_3}(\Xi_1,\Xi'_2,\Xi_3) \times S_{p_1,p_2,p_3-1}(\Xi_1,\Xi_2,\Xi'_3), \nonumber \\
\Xhath{2} \simeq & S_{\bf p}\left(\boldsymbol \Xi; (2,3)\right) \times S_{\bf p}\left(\boldsymbol \Xi; (1,3)\right) \times S_{\bf p}\left(\boldsymbol \Xi; (1,2)\right) \nonumber \\
= & S_{p_1,p_2-1,p_3-1}(\Xi_1,\Xi'_2,\Xi'_3) \times S_{p_1-1,p_2,p_3-1}(\Xi'_1,\Xi_2,\Xi'_3) \times S_{p_1-1,p_2-1,p_3}(\Xi'_1,\Xi'_2,\Xi_3), \nonumber \\
\Xhath{3} \simeq & S_{\bf p}\left(\boldsymbol \Xi; (1, 2, 3) \right) \nonumber 
= S_{p_1-1,p_2-1,p_3-1}(\Xi'_1,\Xi'_2,\Xi'_3), \nonumber
\end{align}
and the exterior derivatives coincide with the standard gradient, curl and divergence operators.

For the solution of \eqref{eq:continuousproblem} and \eqref{eq:evp_with_forms} it will be necessary to compute the exterior derivative of spline differential forms. The discrete de Rham sequence guarantees existence of a matrix representation of the exterior derivative $d^k$, for $k = 0, \ldots, n-1$.
Following \cite{HIP99c} we denote by $\mathbf{D}^k$ these rectangular matrices\footnote{In the three-dimensional case, they are usually denoted by $\mathbf{G}, \mathbf{C}$ (or $\mathbf{R}$) and $\mathbf{D}$, see for instance \cite[Chapter~5]{BOS98b}.}. Thanks to the choice of the basis functions they are sparse incidence matrices, with only $1$ and $-1$ nonzero entries, defined on an auxiliary Cartesian mesh, also called the Greville mesh \cite{ratnani,BSV14,evans_hierarchical_2020}.
We remark that this auxiliary mesh does not need to be constructed in practice. With some abuse of notation, we will also denote by $\mathbf{D}^k$ the matrices when restricting the operators to subspaces $X^k_{h,0}$.


\section{The discrete Hodge--dual complex} \label{sec:dsplinecomplex}

As already mentioned in the introduction, we look for a method based on two de Rham complexes of splines. 
This section deals therefore with constructing the dual sequence and contextually how to connect the sequences to make the final system solvable.

\subsection{The dual spline complex}
Relying on the high continuity of splines, we will define a dual complex without introducing a dual mesh. We will first explain the construction in the one-dimensional case, and then extend it to arbitrary dimensions (in principle) by tensorization. The idea of construction of a dual complex with spline spaces appeared first in \cite{hiemstra_isogeometric_2011} and was also suggested in \cite{BSV14}. It was then used in \cite{Buffa_2020aa,kapidani2021treecotree} to define a Lagrange multiplier for mortar gluing.

\subsubsection{The one-dimensional dual complex}
From here on we assume that the degree is greater than one and that the 0--forms of the primal complex are at least $C^1$ continuous.
We start on the unit interval, where the exterior derivative amounts to the full derivative and the maximum order non-trivial differential forms we can define are therefore 1--forms. For the primal complex, we have the discrete spaces of 0--forms and 1--forms given by
\[
\Xhath{0} = S_{p}(\Xi), \quad \Xhath{1} = \{\hat \omega^h dx : \hat \omega^h \in S_{p-1}(\Xi')\},
\]
and the one-dimensional primal sequence
\begin{equation*}
    \begin{CD}
    \Xhathb{0} @>d^0>> \Xhathb{1},
    \end{CD}
    \label{eq:cd-iga-1d}
\end{equation*}
where the zero subscript denotes homogeneous boundary conditions as in \eqref{eq:cd-iga-tp-bc}.

In principle, functions in the primal space of discrete 1--forms $\Xhath{1}$ are only required to be in the $L^2$ space. But since we assume that functions in $\Xhath{0}$ are $C^1$ continuous, the functions in $\Xhath{1}$ are at least continuous. We can thus construct a new sequence starting in which the spline space $S_{p-1}(\Xi')$ is the space of 0--forms. We define
\[
\XhathT{0} = S_{p-1}(\Xi'), \quad \XhathT{1} = \{\hat \omega^h dx : \hat \omega^h \in S_{p-2}(\Xi'')\},
\]
and we have the one-dimensional dual sequence
\begin{equation*}
    \begin{CD}
    \XhathT{0} @>d^0>> \XhathT{1},
    \end{CD}
    \label{eq:cd-iga-1d-dual}
\end{equation*}
where now vanishing boundary conditions are not included anymore. In practice, the definition of this dual complex is not different from the one we have already seen in \eqref{eq:cd-iga-tp-bc}.
By construction, there is a trivial isomorphism between $\Xhathb{1}$ and $\XhathT{0}$. Indeed, in the unit interval their corresponding spaces of proxy fields are identical, as they are both $S_{p-1}(\Xi')$. Moreover, due to the vanishing boundary conditions, the algebraic dimension of the spaces $\Xhathb{0}$ and $\XhathT{1}$  also matches making the definition of an isomorphism in the finite dimensional setting possible. However, its optimal explicit construction in terms of a discrete Hodge--star operator, which will occupy Section~\ref{sec:discrete_hodge}, is not immediate. Note also that if there is a nontrivial metric on the physical domain, i.e. on a one-dimensional manifold, the push-forward makes also the spaces $X_{h,0}^1$ and $\XhT{0}$, defined in the physical domain, different.

\subsubsection{The tensor-product dual complex}
The construction for dimension $n>1$ proceeds by tensorization of the one-dimensional case. As before, we first define the primal space of 0--forms as $\Xhath{0} = S_{\mathbf{p}}(\boldsymbol{\Xi})$, which gives the primal sequence as in \eqref{eq:cd-iga-tp-bc},
\begin{equation*}
    \begin{CD}
   \Xhathb{0}   @>d^0>>  \Xhathb{1}  @>d^1>> \dots @>d^{n-2}>>  \Xhathb{n-1} @>d^{n-1}>> \Xhathb{n}.
    \end{CD}
\end{equation*}
Assuming that the spline functions in $\Xhath{0}$ are $C^1$ continuous, the functions defining the space of $n$--forms are at least continuous. Hence, as in the one dimensional case, we can use the same spline space appearing in primal $n$--forms to define the dual space of 0--forms, and to start a second sequence. Therefore, with obvious notation, we define the dual space of 0--forms as 
\[
\XhathT{0} = \otimes_{\ell=1}^n S_{p_\ell-1}(\Xi'_\ell) = S_{\mathbf{p-1}}(\boldsymbol{\Xi}'),
\]
from which we have the dual complex
\begin{equation*}
    \begin{CD}
   \XhathT{0}   @>d^0>>  \XhathT{1}  @>d^1>> \dots @>d^{n-2}>>  \XhathT{n-1} @>d^{n-1}>> \XhathT{n}. \\
    \end{CD}
\end{equation*}
The same arguments that we used in the one-dimensional case apply via tensorization: each parametric direction for which the univariate splines in the definition $\Xhathb{k}$ have degree $p$ corresponds to the parametric direction for which univariate splines in  $\XhathT{n-k}$ have degree $p-2$, while all remaining directions are of equal degree $p-1$. 
This is easily seen again taking the vector proxy field example in three dimensions, lowering all degrees by one and comparing the appropriate $(n-k,k)$ pairs. Consequently, the space $\Xhathb{k}$ and its ``dual'' $\XhathT{n-k}$ have the same dimension.

In practice, the construction of this dual complex does not differ from the one already explained in Section~\ref{sec:psplinecomplex}, as we are simply defining a spline sequence of one degree lower than the primal one, using the same knot vector. Consequently, the explanations therein are also valid for the dual complex. In particular, we can map discrete spaces of differential forms to the physical domain $\Omega$, to define the spaces $\XhT{k}$ which form the dual complex in $\Omega$. Moreover, the exterior derivative across all the dual complex is well defined, and given by incidence matrices on another auxiliary control mesh. To distinguish them from the operators of the primal complex, we denote them with $\widetilde{\mathbf{D}}^k$.

\begin{remark}
All spaces in the two sequences are defined on the same mesh (given by the knot vector), and there is no need to explicitly construct a dual mesh. 
This achievement is thanks to the high continuity of the spline spaces in the primal complex, and it is not reproducible in standard FEM approaches.
It is also worth to remark that the construction works for any degree, any regularity, and for non-uniform knot vectors, as long as the starting space in the primal sequence is $C^1$ continuous.
\end{remark}

\subsubsection{Pairing matrices}
Following \cite{HIP99c}, we can define pairing matrices between the spaces of primal differential $k$--forms $X^k_{h,0}$ and dual differential $(n-k)$--forms $\XhT{n-k}$. 
Let $\omega_h^k \in X^k_{h,0}$ and $\eta_h^{n-k} \in \XhT{n-k}$ be respectively represented by the vectors of degrees of freedom $\boldsymbol{\omega}$ and $\boldsymbol{\eta}$. The square pairing matrices ${\bf K}_k^{n-k}$ and $\tilde{\bf K}^k_{n-k}$, are respectively determined by
\begin{align*}
\boldsymbol{\eta}^\top {\bf K}_k^{n-k} \boldsymbol{\omega} = \int_\Omega \eta_h^{n-k} \wedge \omega_h^k,  \quad \text{ for all } \omega_h^k \in X^k_{h,0}, \eta_h^{n-k} \in \XhT{n-k} , \\
\boldsymbol{\omega}^\top \tilde{\bf K}^k_{n-k} \boldsymbol{\eta} = \int_\Omega \omega_h^{k} \wedge \eta_h^{n-k}, \quad \text{ for all } \omega_h^k \in X^k_{h,0}, \eta_h^{n-k} \in \XhT{n-k}, 
\end{align*}
and from the properties of the wedge product, a fundamental property linking them (see~\cite{HIP99c}) immediately follows:
\begin{equation} \label{eq:pairing_transpose}
    \tilde{\bf K}^k_{n-k} = (-1)^{k (n-k)} ({\bf K}_k^{n-k})^\top.
\end{equation}
Moreover, because the pairing $(X^k_{h,0}, \XhT{n-k})$ is stable, as proved in the framework of mortar methods in \cite{Buffa_2020aa}, all the pairing matrices are invertible.

There is an intimate connection between pairing matrices and discrete exterior derivatives. Recall that we assume vanishing boundary conditions for all spaces of forms in the primal sequence, then from the definition of the pairing matrices and the incidence matrices representing the exterior derivatives, it holds   
\begin{equation*}
    ({\bf D}^{k-1})^\top \tilde{\bf K}^k_{n-k} = (-1)^{k} \tilde{\bf K}^{k-1}_{n-k+1} \tilde{\bf D}^{n-k},
\end{equation*}
for any $k$ and conversely, using \eqref{eq:pairing_transpose} it holds
\begin{equation*}
    (\tilde{\bf D}^{n-k})^\top {\bf K}^{n-k+1}_{k-1} = (-1)^{n-k+1} {\bf K}^{n-k}_{k} {\bf D}^{k-1},
\end{equation*}
which both represent the discrete counterpart of the integration by parts formula in \eqref{eq:int_by_parts}. We remark that as such, no metric is yet involved in the above properties.

\subsection{Discrete Hodge--star operators} \label{sec:discrete_hodge}
We will hereafter introduce three different choices of the discrete Hodge--star operator which preserve the high order nature of the spline based discretization. The first two are in some sense natural.
They are induced by the properties of the Hodge--star operator in Section~\ref{sec:derham-hodge}, and adapt to our setting the discrete Hodge operators as presented in \cite{HIP99c}. Their main drawback is that they are global operators. The third choice we present is instead based on local projection operations from \cite{LLM01}.

\subsubsection{Global discrete Hodge operators}
Let $\omega_h^k \in X^k_{h,0}$, and let us denote in general $u_h^{n-k} = \star_{{h,\gamma}}^k \omega_h^k$, with discrete Hodge--star operator $\star_{{h,\gamma}}^k:X^k_{h,0} \rightarrow \XhT{n-k}$ yet undefined. The first discrete operator, mimicking the definition of the continuous Hodge--star operator in \eqref{eq:def_hodge_diff}, is uniquely determined by
\begin{equation} \label{eq:def_hodge_discrete1}
    \int_\Omega \eta_h^k \wedge \star_{h,\gamma}^k \omega_h^k = \left( \eta_h^k, \omega_h^k \right)_{L_\gamma^2 \Lambda^k(\Omega)} ,
    \quad \forall \eta^k \in X_{h,0}^k,
\end{equation}
\noindent which in matrix form, using boldface symbols for vectors of degrees of freedom, reads
$$\boldsymbol{\eta}^\top \tilde{\bf K}^k_{n-k} {\bf u} = \boldsymbol{\eta}^\top  {\bf M}^k_\gamma \boldsymbol \omega.$$
\noindent where ${\bf M}^k_\gamma$ is the standard mass matrix for $X^k_{h,0}$ and $\tilde{\bf K}^k_{n-k}$ is the pairing matrix between $X^k_{h,0}$ and $\XhT{n-k}$, which we recall being square and invertible. Our first discrete Hodge--star operator is hence written in matrix form as 
\begin{equation*}
    {\bf H}_\gamma^k = (\tilde{\bf K}^k_{n-k})^{-1}  {\bf M}^k_\gamma,
\end{equation*}
\noindent whereas, for the definition of the second discrete Hodge--star operator, we mimic property \eqref{eq:def_hodge_diff_inv} from the continuous setting: the operator is this time uniquely determined by
\begin{equation}\label{eq:def_hodge_discrete2}
(\eta_h^{n-k},\star_{h,\gamma}^k \omega_h^k)_{L^2_{1/\gamma}\Lambda^{n-k}(\Omega)} = (-1)^{k(n-k)} \int_\Omega \eta_h^{n-k} \wedge \omega_h^k, \quad \forall \eta_h^{n-k} \in \XhT{n-k},
\end{equation}
\noindent from which it ensues, in matrix form again,
$$\boldsymbol{\eta}^\top\tilde{\bf M}^{n-k}_{1/\gamma} {\bf u} = (-1)^{k(n-k)} \boldsymbol{\eta}^\top {\bf K}_k^{n-k} \boldsymbol \omega,$$
\noindent from which, in turn, we deduce its matrix representation as 
\begin{equation*}
    {\bf H}_\gamma^k = (-1)^{k(n-k)} (\tilde {\bf M}^{n-k}_{1/\gamma})^{-1} {\bf K}_k^{n-k},
\end{equation*}
\noindent where, similarly to above $\tilde {\bf M}^{n-k}_{1/\gamma}$ is the mass matrix for the space $\XhT{n-k}$ and ${\bf K}_k^{n-k}$ is the pairing matrix.

\subsubsection{A local discrete Hodge operator} \label{sec:LLM}
The global nature of the two above choices is somehow unfortunate since methods based on dual grids use sparse, sometimes even diagonal, Hodge operators. We propose a third version of the discrete Hodge operator which partially mends this drawback.
By glancing at the second proposed choice for the operator, it is easy to see that
\[
\star_{{h,\gamma}}^k \omega_h^k = \tilde \Pi^{n-k}_{L^2} (\star^k_\gamma \omega^k_h),
\]
i.e., the discrete Hodge operator is a standard $L^2$--projection composed with the continuous Hodge--star operator.

The idea is then to use another projector, which has approximation properties proven to be equivalent to the $L^2$--projection, but exhibiting a local nature. We here use a quasi-interpolant projector for splines, which we label $\tilde \Pi^{n-k}_{\text{\tiny LLM}}$, which enters into the family of quasi-interpolants defined in \cite{LLM01}.
The projector is best described as an algorithm in three steps, which for $\star_{{h,\gamma}}^k$ acts on the basis of the target space $\XhT{n-k}$ as follows:

\begin{enumerate}
\item[i)] For each basis function $\lambda_{i}$ of the space $\XhT{n-k}$, choose $Q_i \subset \mathrm{supp} \lambda_{i}$.
\item[ii)] Project locally on $Q_i$, more precisely on $\XhT{n-k}|_{Q_i}$ (in our case again with an $L^2$--projection $\tilde \Pi^{Q_i}_{L^2} (\star^k_\gamma \omega^k_h|_{Q_i})$).
\item[iii)] Choose the coefficients of the local projection associated to $\lambda_{i}$. This will be the coefficient in the final matrix representation of the global projection.
\end{enumerate}

Steps i) and ii) entail some liberty in the choice of a specific subset of the support of basis functions and on the choice of projection, respectively. We choose local $L^2$--projections in step ii) instead of point interpolation used for the examples in \cite{LLM01}, because it suits better the approximation of general $k$--forms and the definition of the Hodge--star operator, as point evaluation is only meaningful for $0$--forms. For the choice of $Q_i$ in step i), we choose the barycentric element in the support of the univariate B-spline for each specific parametric direction. The $n$-dimensional support $Q_i$ for the local projection will then be the Cartesian product of the one-dimensional supports, mapped to the physical domain, where the projection is finally computed taking into account the Jacobian of the mapping of $\Omega$ (differently from \cite{LLM01}, where the projection is performed on the parametric domain).
We try to take the smallest possible support centered around the maximum of each target basis spline.
If the number of mesh elements in a specific direction is even, e.g. if the degree of the splines is odd or we have repeated knots, we choose the two elements containing the barycentric knot of the support, which tends to reduce the sparsity of the operator.
A matrix representation of the discrete Hodge--star operator is also possible in this case, although an explicit expression is not as simple as for the previous two operators. However, the advantage is that since the operator is local, the resulting matrix is expected to be sparser.

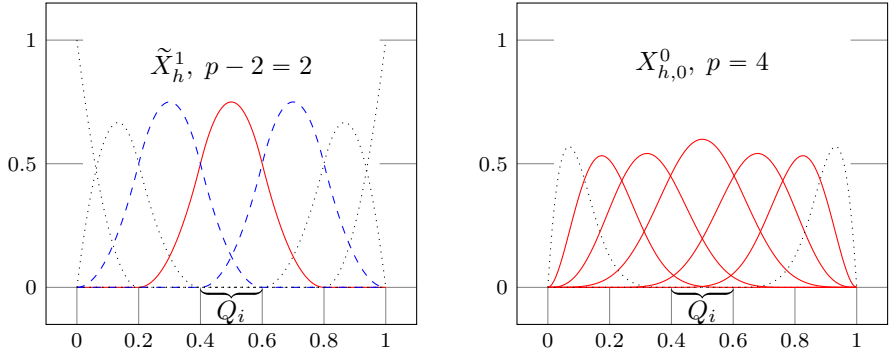
\begin{figure}
    \centering
    \begin{tikzpicture}[]
        \begin{axis}[width=0.5\columnwidth, 
            major tick length = 14pt,
            xtick = {                   0,
            0.200000000000000,
            0.400000000000000,
            0.600000000000000,
            0.800000000000000,
            1.000000000000000}, 
            xtick pos = bottom,
            ytick = {0,0.5,1.0},
            height=0.45\columnwidth, ylabel={},xlabel={},  
            ymin=-0.15, ymax=1.15, ylabel near ticks, tick label style={font=\footnotesize},  legend style={font=\footnotesize, }, label style={font=\footnotesize}, every x tick scale label/.style={at={(1,0)}, anchor=north, yshift=-5pt,inner sep=0pt},legend pos= south east]
        \addplot [red, domain=0:1, range=-0.2:1] table [x index=0, y index=4, col sep=comma] {figures/bsplines_p2r1.csv};
        \addplot [dotted, black, domain=0:1] table [x index=0, y index=1, col sep=comma] {figures/bsplines_p2r1.csv};
        \addplot [dotted, black, domain=0:1] table [x index=0, y index=2, col sep=comma] {figures/bsplines_p2r1.csv};
        \addplot [dashed, blue, domain=0:1] table [x index=0, y index=3, col sep=comma] {figures/bsplines_p2r1.csv};
        \addplot [dashed, blue, domain=0:1] table [x index=0, y index=5, col sep=comma] {figures/bsplines_p2r1.csv};
        \addplot [dotted, black, domain=0:1] table [x index=0, y index=6, col sep=comma] {figures/bsplines_p2r1.csv};
        \addplot [dotted, black, domain=0:1] table [x index=0, y index=7, col sep=comma] {figures/bsplines_p2r1.csv};
        \draw[red,decorate,decoration={calligraphic brace,amplitude=3pt},thick] (axis cs:0.6,-0.01) -- (axis cs:0.4,-0.01);
        \node[] at (axis cs:0.5, 0.9)   (a) {$\tilde{X}_h^1,\;p-2=2$};
        \node[] at (axis cs:0.5, -0.1)   (q) {$Q_i$};
        \end{axis}
    \end{tikzpicture}
\qquad
    \begin{tikzpicture}[]
    \begin{axis}[width=0.5\columnwidth, 
        major tick length = 14pt,
        xtick = {                   0,
        0.200000000000000,
        0.400000000000000,
        0.600000000000000,
        0.800000000000000,
        1.000000000000000},
        xtick pos = bottom,
        ytick = {0,0.5,1.0},
        height=0.45\columnwidth, ylabel={},xlabel={}, xrange=0:1, 
        ymin=-0.15, ymax=1.15, ylabel near ticks,tick label style={font=\footnotesize},  legend style={font=\footnotesize, }, label style={font=\footnotesize}, every x tick scale label/.style={at={(1,0)}, anchor=north, yshift=-5pt,inner sep=0pt},legend pos= south east]
    \addplot [dotted, black, domain=0:1] table [x index=0, y index=2, col sep=comma] {figures/bsplines_p4r3.csv};
    \addplot [red, domain=0:1] table [x index=0, y index=3, col sep=comma] {figures/bsplines_p4r3.csv};
    \addplot [red, domain=0:1] table [x index=0, y index=4, col sep=comma] {figures/bsplines_p4r3.csv};
    \addplot [red, domain=0:1] table [x index=0, y index=5, col sep=comma] {figures/bsplines_p4r3.csv};
    \addplot [red, domain=0:1] table [x index=0, y index=6, col sep=comma] {figures/bsplines_p4r3.csv};
    \addplot [red, domain=0:1] table [x index=0, y index=7, col sep=comma] {figures/bsplines_p4r3.csv};
    \addplot [dotted, black, domain=0:1] table [x index=0, y index=8, col sep=comma] {figures/bsplines_p4r3.csv};
    \draw[red,decorate,decoration={calligraphic brace,amplitude=3pt},thick] (axis cs:0.6,-0.01) -- (axis cs:0.4,-0.01);
    \node[] at (axis cs:0.5, 0.9)   (a) {${X}_{h,0}^0,\;p=4$};
    \node[] at (axis cs:0.5, -0.1)   (q) {$Q_i$};
    \end{axis}
    \end{tikzpicture}
\caption{Choice of the support of integration for the local projector in the case of a discrete Hodge--star operator $\star_h^0 : {X}_{h,0}^0\mapsto\XhT{1}$, from splines of degree four to splines of degree two. For the solid red B-spline 1--form $\lambda_i \in \XhT{1}$ highlighted in the left panel, $Q_i$ is the central element. Blue dashed B-splines contribute to the matrix involved in the local $L^2$ projection operation. Solid red B-splines in the starting space $X_{h,0}^0$ on the right panel are the basis functions that do not vanish in $Q_i$, and may have a nonzero contribution from $\lambda_i$ after applying the Hodge--star operator. The rest of the basis functions in both panels are shown dotted in black.}
    \label{fig:bsplines}
\end{figure} 

To better understand the sparsity of this local Hodge operator, let us consider the one-dimensional case with the maximum allowed smoothness, and assume we construct a Hodge--star operator $\star_{h,\gamma}^0: X_{h,0}^0\mapsto \XhT{1}$, from the space of primal $0$--forms of degree $p$ to the space of dual $1$--forms of degree $p-2$. Each basis function $\lambda_i$ of $S_{p-2}(\Xi'')$ is supported on $p-1$ knot intervals, therefore we choose $Q_i$ as the central knot interval if $p$ is even, see Fig.~\ref{fig:bsplines}, or as the union of two central knot intervals if $p$ is odd. Near the endpoints of the knot vector the choice has to be adapted due to the reduced support of the basis functions. 
Considering that each element is contained in the support of $p+1$ functions of $X_{h,0}^0$, and $p-1$ functions of $\XhT{1}$, in the worst case of choosing two elements each local $L^2$ projection requires solving a $p \times p$ system of equations, with a multi-column right-hand side of size at most $p+2$ columns. Each system can be solved concurrently, and the resulting coefficients associated with $\lambda_i$ will contribute to $p+2$ nonzero entries in each row of the final matrix representation of the projector. In higher dimensional cases, the accurate computation of the sparsity is given by tensor product arguments on univariate spaces, noting that for $k$--forms with $0<k<n$, the presence of a non-diagonal Jacobian matrix will reduce sparsity by mixing up forms with different multi-index $\boldsymbol{\alpha}$ (or different components of vector proxy fields).

\subsection{The discrete method}
The discrete version of the continuous problem \eqref{eq:continuousproblem} is given by the problem of finding $\omega_h \in X^k_{h,0}$ such that 
\begin{equation} \label{eq:discreteproblem}
(-1)^{k+1} d^{n-k-1} \star_{h,\gamma}^{k+1} (d^k \omega_h) + \star^k_{h,\beta} \omega_h = \star^k_h f,
\end{equation}
where the discrete Hodge--star operators can follow any of the three constructions given above. Using the matrix representation of the exterior derivatives and the discrete Hodge--star operators, the problem can also be written in matrix form as
\begin{equation}  \label{eq:problem_matrix}
(-1)^{k+1} \tilde{{\bf D}}^{n-k-1} {\bf H}_\gamma^{k+1} ({\bf D}^k \boldsymbol{\omega}_h) + {\bf H}^k_\beta \boldsymbol{\omega}_h = {\bf H}^k \boldsymbol{f}_h,
\end{equation}
where we remark that the resulting system matrix is square even though the approximation is not performed in Galerkin fashion, due to the dimension of the spaces of the dual complex. With a completely analogous procedure, we obtain the discrete version of the eigenvalue problem \eqref{eq:evp_with_forms}.

Finally, we note that the right-hand side in \eqref{eq:discreteproblem} must be a discrete differential form in $\XhT{n-k}$. Although we have written the equation applying the same kind of discrete Hodge--star operator, nothing prevents us to use, for instance, an $L^2$ projection into $\XhT{n-k}$.
Both choices provide a consistent right-hand side.

\section{Error analysis of the elliptic problem}\label{sec:erroranalysis}

In this section we perform the error analysis of the solution of the discrete problem \eqref{eq:discreteproblem}. Unlike \cite{HIP99c}, for the analysis of the problem we leverage standard variational techniques. We will see that the first Hodge operator is in fact equivalent to a standard Galerkin discretization, while the other two yield new methods with a reduced order of convergence. For the sake of simplicity, in this section we will consider same starting degree in every direction, i.e., $p_\ell = p$ for $\ell = 1, \ldots, n$.

We first multiply (using the wedge product) the equation of the continuous problem \eqref{eq:continuousproblem} by a test function, and integrate on $\Omega$, to obtain the variational formulation
\begin{equation} \label{eq:problem_continuous}
a(\omega, \eta) = {\cal F}(\eta) \quad \text{ for all } \eta \in H_0 \Lambda^k (\Omega),
\end{equation}
where, using \eqref{eq:int_by_parts} and \eqref{eq:def_hodge_diff}, we get three equivalent expressions for the bilinear form $a(\cdot, \cdot)$, namely:
\begin{align} \label{eq:definition_a}
a(\omega, \eta) := & (-1)^{(k+1)} \int_\Omega \eta \wedge d^{n-k-1} \star^{k+1}_\gamma (d^k \omega)  + \int_\Omega  \eta \wedge \star^k_\beta \omega \notag \\
 = & 
 \int_\Omega  d^k \eta \wedge \star^{k+1}_\gamma (d^k \omega) + \int_\Omega \eta \wedge \star^k_\beta \omega \\
 = & \left(d^k \omega, d^k \eta \right)_{L^2_\gamma \Lambda^{k+1}(\Omega)} + \Big( \omega, \eta \Big)_{L^2_\beta\Lambda^k(\Omega)}, \notag
\end{align}
while the linear form ${\cal F}$ is defined as
\[
{\cal F}(\eta) := \int_\Omega \eta \wedge \star^k f = (f, \eta)_{L^2\Lambda^k(\Omega)}.
\]
The discrete Hodge operators define analogous bilinear forms and linear forms.
We denote by $a_h$ and ${\cal F}_h$ for all of them. Since the exterior derivative is well defined for spline spaces of the two sequences, integration by parts remains valid and they are defined by
\begin{align}
 a_h(\omega_h, \eta_h) := & (-1)^{(k+1)} \int_\Omega \eta_h \wedge d^{n-k-1} \star^{k+1}_{h,\gamma} (d^k \omega_h) + \int_\Omega \eta_h \wedge \star^k_{h,\beta} \omega_h \notag \\
 = & \int_\Omega d^k \eta_h \wedge \star^{k+1}_{h,\gamma} (d^k \omega_h) + \int_\Omega \eta_h \wedge \star^k_{h,\beta} \omega_h  \label{eq:definition_ah} \\
 {\cal F}_h(\eta_h) := & \int_\Omega \eta_h \wedge \star^k_h f , \notag
\end{align}
and the discrete problem \eqref{eq:discreteproblem} is then equivalent to find $\omega_h \in X^k_{h,0}$ such that
\begin{equation} \label{eq:discreteproblem_weak}
a_h(\omega_h, \eta_h) = {\cal F}_h (\eta_h) \quad \forall \eta_h \in X^k_{h,0}.
\end{equation}

Clearly, the bilinear form $a(\cdot,\cdot)$ is symmetric, continuous and coercive in $H_0\Lambda^k(\Omega)$, and the linear form ${\cal F}$ is continuous. The analysis is then based on first Strang's lemma \cite[Thm.~5.5.1]{Quarteroni-Valli-1994}, that we restate here for convenience.
\begin{lemma}
Assuming that the bilinear form $a_h(\cdot, \cdot)$ is uniformly coercive over $X^k_{h,0}$, then there exists $\omega_h \in X^k_{h,0}$ a unique solution to \eqref{eq:discreteproblem_weak}, and it holds that
\begin{align*}
\| \omega - \omega_h \|_{H\Lambda^k(\Omega)} & \le 
C_1 \sup_{\lambda_h \in X^k_{h,0} \setminus\{0\}} \frac{|{\cal F}(\lambda_h) - {\cal F}_h(\lambda_h)|}{\|\lambda_h\|_{H\Lambda^k(\Omega)}}\\
& +\inf_{\eta_h \in X^k_{h,0}}\!\left[ C_2 \| \omega - \eta_h \|_{H\Lambda^k(\Omega)} + C_3\! \sup_{\lambda_h \in X^k_{h,0} \setminus\{0\}} \!\hspace{-6pt}\frac{|a(\eta_h,\lambda_h) - a_h(\eta_h,\lambda_h)|}{\|\lambda_h\|_{H\Lambda^k(\Omega)}} \right],
\end{align*}
where constants $C_1, C_2$ and $C_3$ depend on the coerciveness constant of $a_h(\cdot,\cdot)$, and $C_2$ also depends on the continuity constant of $a(\cdot,\cdot)$.
\end{lemma}

In all the three cases, that is, independently of the discrete Hodge--star operator, the second term of the error can be bounded using the projectors in \cite{BRSV11}, and for $0 \le s \le p$ we have
\begin{equation} \label{eq:interpolation_error}
\inf_{\eta_h \in X^k_{h,0}} \| \omega - \eta_h \|_{H\Lambda^k(\Omega)} \le C h^{\min\{p,s\}} \| \omega \|_{H^s\Lambda^k(\Omega)}, \quad \text{ for } \omega \in H_0\Lambda^k(\Omega) \cap H^s\Lambda^k(\Omega).
\end{equation}
Therefore, the error analysis reduces to the study of the consistency errors. From now on, and to simplify the analysis, we will assume that there is no consistency error in the right-hand side, i.e. ${\cal F}(\eta_h) = {\cal F}_h(\eta_h)$. 
This is a reasonable hypothesis, since $f$ is a given datum, and not a discrete function.
\begin{remark}
The correct choice of the right-hand side ${\cal F}_h(\eta_h)$, or equivalently of $f_h \equiv \star^k_h f$, may be relevant in some cases, for instance when the discrete function $f_h$ must be divergence free.
\end{remark}

The first theoretical result concerns the first Hodge-star operator, defined in \eqref{eq:def_hodge_discrete1}. In fact, for this operator there is no consistency error, and the total error behaves as the interpolation error in \eqref{eq:interpolation_error}. Therefore, as long as there is also no consistency error in the right-hand side, the solutions obtained applying this first operator or a standard Galerkin method are equivalent, although yielding different matrices.
\begin{theorem} \label{thm:first-hodge}
Let $\omega \in H_0 \Lambda^k(\Omega) \cap H^s \Lambda^k(\Omega)$ be the solution to \eqref{eq:problem_continuous}. Then, the solution $\omega_h$ of the discrete problem \eqref{eq:discreteproblem_weak}, with the discrete Hodge--star operator defined in \eqref{eq:def_hodge_discrete1}, satisfies
\[
\|\omega - \omega_h \|_{H\Lambda^k(\Omega)} \le C h^{\min\{s,p\}} \| \omega \|_{H^s\Lambda^k(\Omega)}.
\]
\end{theorem}
\begin{proof}
Using the definition of the discrete Hodge--star operator \eqref{eq:def_hodge_discrete1} in \eqref{eq:definition_ah}, it holds that
\[
a_h(\eta_h, \lambda_h) = a(\eta_h, \lambda_h) \text{ for all }\eta_h, \lambda_h \in X^k_{h,0},
\]
and the result immediately follows from first Strang's lemma.
\end{proof}

For the analysis of the second and third Hodge--star operators, we first write them as the combination of the continuous Hodge--star operator with a projection onto $\XhT{n-k}$, which we can denote in general by $\tilde{\Pi}^{n-k}$, resulting in
\begin{equation} \label{eq:hodge-projector}
\star^k_{h,\beta} \omega_h = \tilde{\Pi}^{n-k} (\star^k_\beta \omega_h),
\end{equation}
and analogously for $\star^{k+1}_{h,\gamma}$. Both operators have the approximation properties stated in the following lemma, which we present without proof. For the $L^2$-projection the result is due to the approximation properties of B-splines, while for the local operator it comes from the fact that this is nevertheless based on local $L^2$ projections, possibly with a weight due to the presence of the transformation to the physical domain.
\begin{lemma} \label{lemma:approx}
Let $\tilde \Pi^{n-k}$ denote either $\tilde \Pi^{n-k}_{L^2}$, the $L^2$-projection onto $\XhT{n-k}$, or the local projection $\tilde \Pi^{n-k}_{\text{\tiny LLM}}$ defined in Section~\ref{sec:LLM}. Then, for any $\eta \in H^s\Lambda^{n-k}(\Omega)$ with $0 \le s \le p$ it holds that
\begin{align*}
& \| \eta - \tilde \Pi^{n-k} (\eta) \|_{L^2\Lambda^{n-k}(\Omega)} \lesssim h^{\min\{s,p-1\}} \| \eta \|_{H^s \Lambda^{n-k}(\Omega)}.
\end{align*}
\end{lemma}
Note that the order of the approximation is reduced with respect to \eqref{eq:interpolation_error}, because we are projecting onto spaces of a spline sequence of one degree less. 
Before analyzing the error, we shall prove that the bilinear forms defined by these discrete Hodge--star operators are coercive.
\begin{lemma} \label{lem:coerciveness}
There exists $h_0 >0$ such that, for $h < h_0$, 
\[
a_h(\eta_h, \eta_h) \gtrsim \| \eta_h \|_{H \Lambda^{k}(\Omega)} \quad \text{ for all } \eta_h \in X^k_{h,0}.
\]
\end{lemma}
\begin{proof}
Using first the definition of $a_h(\cdot, \cdot)$ and expression \eqref{eq:hodge-projector} of the discrete Hodge--star operators, property \eqref{eq:def_hodge_diff_inv} and the fact that $\star^k_\beta \circ \star^{n-k}_{1/\beta} = (-1)^{k(n-k)}$ (and similarly for $\star^{k+1}_\gamma$), we have
\begin{align*}
& a_h(\eta_h, \eta_h) = \int_\Omega d^k \eta_h \wedge \tilde \Pi^{n-k-1} (\star^{k+1}_\gamma (d^k \eta_h)) + \int_\Omega \eta_h \wedge \tilde \Pi^{n-k} (\star^k_\beta \eta_h) \\
& = \left(\star^{k+1}_\gamma d^k \eta_h, \tilde \Pi^{n-k-1} (\star^{k+1}_\gamma (d^k \eta_h))\right)_{L^2_{1/\gamma}\Lambda^{n-k-1}(\Omega)} + 
\left(\star^k_\beta \eta_h, \tilde \Pi^{n-k} (\star^k_\beta \eta_h)\right)_{L^2_{1/\beta}\Lambda^{n-k}(\Omega)} \\
& = \left(\phi_h, \tilde \Pi^{n-k-1} (\phi_h))\right)_{L^2_{1/\gamma}\Lambda^{n-k-1}(\Omega)} + 
\left(\psi_h, \tilde \Pi^{n-k} (\psi_h)\right)_{L^2_{1/\beta}\Lambda^{n-k}(\Omega)},
\end{align*}
where we have introduced $\phi_h := \star^{k+1}_\gamma (d^k \eta_h) \in \XhT{n-k-1}$ and $\psi_h := \star^k_\beta \eta_h \in \XhT{n-k}$, to simplify notation. By adding and subtracting equal terms, applying Cauchy-Schwarz inequality, and finally by virtue of the approximation result of Lemma~\ref{lemma:approx}, which can be used due to all the spline spaces involved in the definition $X^k_{h,0}$ being at least continuous, it ensues
\begin{align*}
a_h(\eta_h, \eta_h) 
& \ge \left(\phi_h, \phi_h\right)_{L^2_{1/\gamma}\Lambda^{n-k-1}(\Omega)} - \left|\left(\phi_h, \phi_h - \tilde \Pi^{n-k-1} \phi_h)\right)_{L^2_{1/\gamma}\Lambda^{n-k-1}(\Omega)}\right| \\
& + \left(\psi_h, \psi_h \right)_{L^2_{1/\beta}\Lambda^{n-k}(\Omega)} - \left|\left(\psi_h, \psi_h - \tilde \Pi^{n-k} \psi_h \right)_{L^2_{1/\beta}\Lambda^{n-k}(\Omega)}\right| \\
& \gtrsim \| \phi_h \|^2_{L^2\Lambda^{n-k-1}(\Omega)} - \| \phi_h\|_{L^2\Lambda^{n-k-1}(\Omega)} \| \phi_h - \tilde \Pi^{n-k-1} \phi_h\|_{L^2\Lambda^{n-k-1}(\Omega)} \\
& + \| \psi_h \|^2_{L^2\Lambda^{n-k}(\Omega)} - \| \psi_h \|_{L^2\Lambda^{n-k}(\Omega)} \| \psi_h - \tilde \Pi^{n-k} \psi_h \|_{L^2\Lambda^{n-k}(\Omega)} \\
& \gtrsim \| \phi_h \|^2_{L^2\Lambda^{n-k-1}(\Omega)} + \| \psi_h \|^2_{L^2\Lambda^{n-k}(\Omega)} \\
& - C h \left( \| \phi_h \|_{L^2\Lambda^{n-k-1}(\Omega)} \| \phi_h \|_{H\Lambda^{n-k-1}(\Omega)} + \| \psi_h \|_{L^2\Lambda^{n-k}(\Omega)} \| \psi_h\|_{H\Lambda^{n-k}(\Omega)} \right),
\end{align*}
where we have also used the equivalence of the standard $L^2$ norm with the weighted norms, due to the assumptions on the material properties $\gamma$ and $\beta$. Coerciveness then holds, for $h$ sufficiently small, by recalling the definitions of $\phi_h$ and $\psi_h$, the definition of the norm $H\Lambda^k(\Omega)$, and by the continuity of the Hodge--star operators and their inverse, using again the assumptions on the material properties.
\end{proof}

\begin{theorem}\label{thm:error_bound}
Let $\omega \in H_0 \Lambda^k(\Omega) \cap H^s \Lambda^k(\Omega)$ be the solution to \eqref{eq:problem_continuous}, and let $\omega_h$ the solution of the discrete problem \eqref{eq:discreteproblem_weak}, with the discrete Hodge--star operators defined in \eqref{eq:def_hodge_discrete2} or in Section~\ref{sec:LLM}. Then, there exists $h_0 > 0$ such that, for $h < h_0$, it holds that
\begin{equation*}
\| \omega - \omega_h \|_{H\Lambda^k(\Omega)} \lesssim h^{\min\{s,p-1\}} \| \omega \|_{H^s\Lambda^k(\Omega)},
\end{equation*}
\end{theorem}
\begin{proof}
The coerciveness of the bilinear form $a_h$ was proved in Lemma~\ref{lem:coerciveness}, for $h$ sufficiently small. The result is then a consequence of Strang's lemma, where we have to analyze the consistency error. As in the proof of Lemma~\ref{lem:coerciveness}, let us introduce $\phi_h := \star^{k+1}_\gamma (d^k \eta_h)$ and $\psi_h := \star^k_\beta \eta_h$ to simplify notation. Using first the definition of the bilinear forms $a(\cdot, \cdot)$ and $a_h(\cdot, \cdot)$ in \eqref{eq:definition_a} and \eqref{eq:definition_ah}, the definition of the discrete Hodge--star operators as the composition of a projector and the continuous Hodge--star operator as in \eqref{eq:hodge-projector}, and subsequently using the property of the continuous Hodge--star operator in~\eqref{eq:def_hodge_diff_inv} (neglecting the signs, as we take absolute values), we get 
\begin{align*}
& |a(\eta_h, \lambda_h) - a_h(\eta_h, \lambda_h)| = 
\left| \int_\Omega  d^k \lambda_h \wedge \left(\phi_h - \tilde{\Pi}^{k+1} (\phi_h)\right) + \int_\Omega \lambda_h \wedge \left(\psi_h - \tilde{\Pi}^{k}(\psi_h)\right) \right| \\
& = \left| (\star^{k+1}_\gamma d^k \lambda_h , \phi_h - \tilde{\Pi}^{k+1} (\phi_h))_{L^2_{1/\gamma}\Lambda^{n-k-1}(\Omega)} + \left( \star^k_\beta \lambda_h , \psi_h - \tilde{\Pi}^{k}(\psi_h) \right)_{L^2_\beta \Lambda^{n-k}(\Omega)} \right|.
\end{align*}
Then, by Cauchy-Schwarz inequality we obtain
\begin{align*}
|a(\eta_h, \lambda_h) - a_h(\eta_h, \lambda_h)| 
\lesssim & \| \star^{k+1}_\gamma d^k \lambda_h \|_{L^2_{1/\gamma}\Lambda^{n-k-1}(\Omega)} \, \| \phi_h - \tilde{\Pi}^{k+1} (\phi_h)\|_{L^2_{1/\gamma}\Lambda^{n-k-1}(\Omega)} \\
+ & \|\star^k_\beta \lambda_h \|_{L^2_{1/\beta} \Lambda^{n-k}(\Omega)} \, \| \psi_h - \tilde{\Pi}^{k}(\psi_h) \|_{L^2_{1/\beta} \Lambda^{n-k}(\Omega)}.
\end{align*}
Recalling the definition of $\phi_h$ and $\psi_h$, the approximation properties of the projection in Lemma~\ref{lemma:approx} plus the continuity of the Hodge--star operators yield
\begin{align*}
& |a(\eta_h, \lambda_h) - a_h(\eta_h, \lambda_h)| 
\lesssim  h^{\min\{s,p-1\}} \| \star^{k+1}_\gamma d^k \lambda_h \|_{L^2_{1/\gamma}\Lambda^{n-k-1}(\Omega)} \, \| \phi_h\|_{H^s \Lambda^{n-k-1}(\Omega)} \\
& \hspace{3.56cm} + h^{\min\{s,p-1\}} \|\star^k_\beta \lambda_h \|_{L^2_{1/\beta} \Lambda^{n-k}(\Omega)} \, \| \psi_h \|_{H^s \Lambda^{n-k}(\Omega)} \\
& \lesssim h^{\min\{s,p-1\}} (\| d^k \lambda_h \|_{L^2_\gamma \Lambda^{k+1}(\Omega)} \| (d^k \eta_h)\|_{H^s \Lambda^{k+1}(\Omega)} + \| \lambda_h \|_{L^2_{\beta} \Lambda^{k}(\Omega)} \| \eta_h \|_{H^s \Lambda^{k}(\Omega)}),
\end{align*}
and the estimate follows from the definition of the norms, the assumptions on the material properties $\gamma$ and $\beta$ given in Section~\ref{sec:derham}, and from using \eqref{eq:interpolation_error} for the other term in Strang's lemma.
\end{proof}

\begin{remark}
For the second and third Hodge--star operators we have only been able to prove coerciveness under the condition of a sufficiently fine mesh. However, all the numerical tests we have conducted so far to estimate the coerciveness constant strongly suggest the bilinear form to be uniformly coercive.
\end{remark}

\section{Numerical results} \label{sec:tests}




In the following the analysis from the previous section is numerically tested on selected elliptic problems, furthermore the important case of eigenvalue approximation is studied numerically to consolidate the more general structure preserving aim of the proposed method. Time dependent problems will be studied in a subsequent companion contribution, for which the present study provides a necessary foundation.

\begin{figure}[!t]
    \centering
    \includegraphics[width=0.5\textwidth]{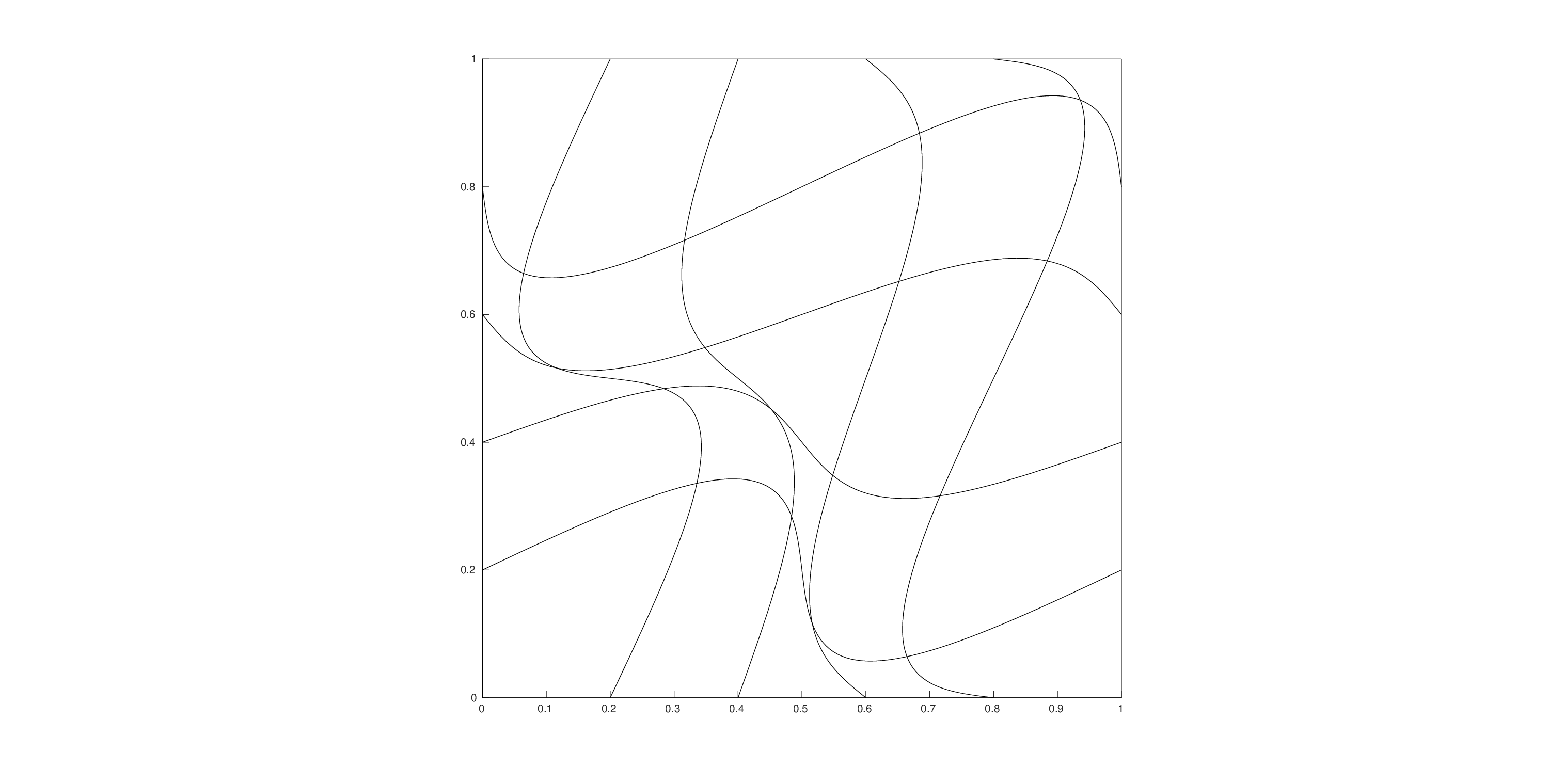}\qquad
    \includegraphics[width=0.2\textwidth,height=0.2\textheight]{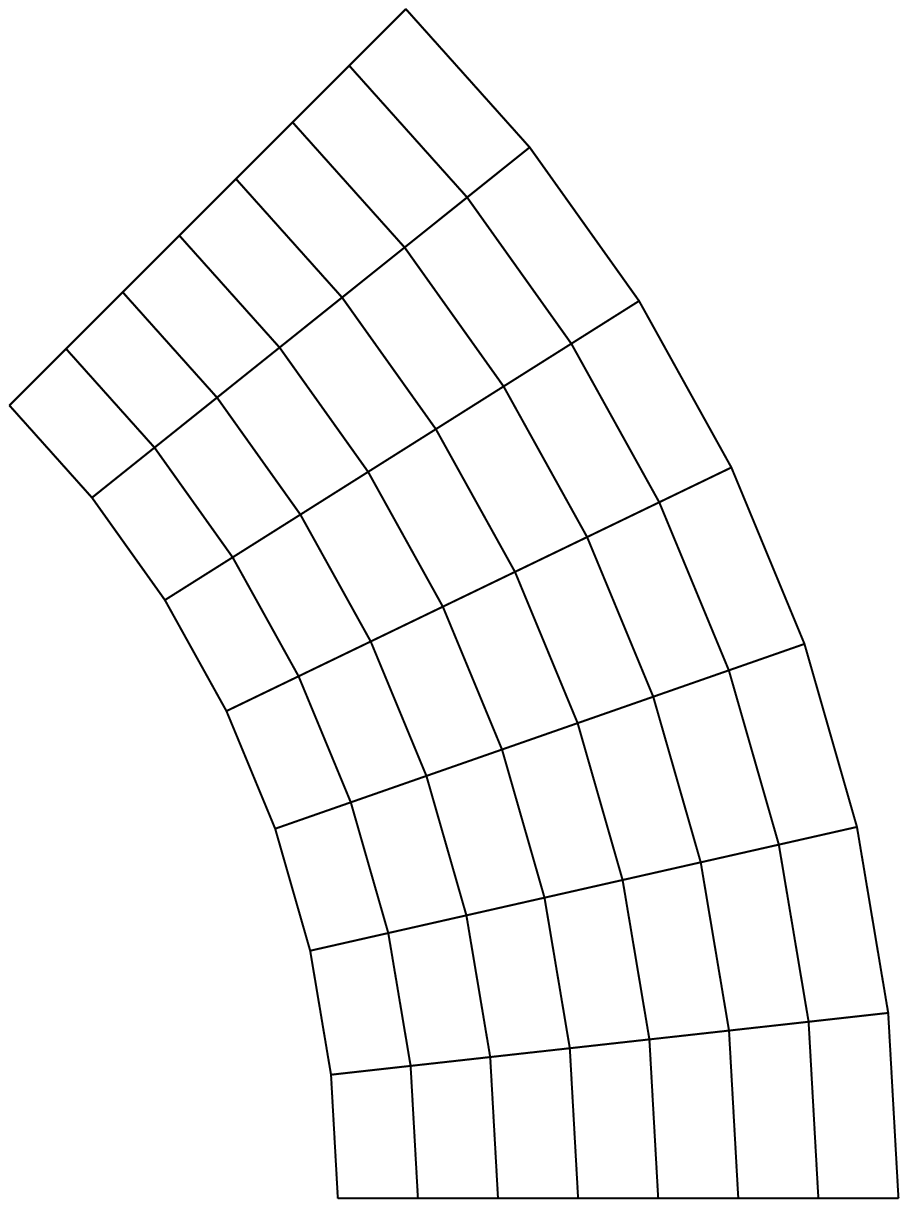}
    \caption{Reparametrization of the unit square via the map in \eqref{eq:harmonic_map} with $a=0.3$ and the one-eighth ring geometry $\Omega = [1,2]\times[0,\pi/4]$ in polar coordinates $(r,\varphi)$, used for both the boundary value and eigenvalue test problem studied.}\label{fig:reparam}
\end{figure}

\subsection{Elliptic boundary value problem} \label{sec:elliptic_numerical}
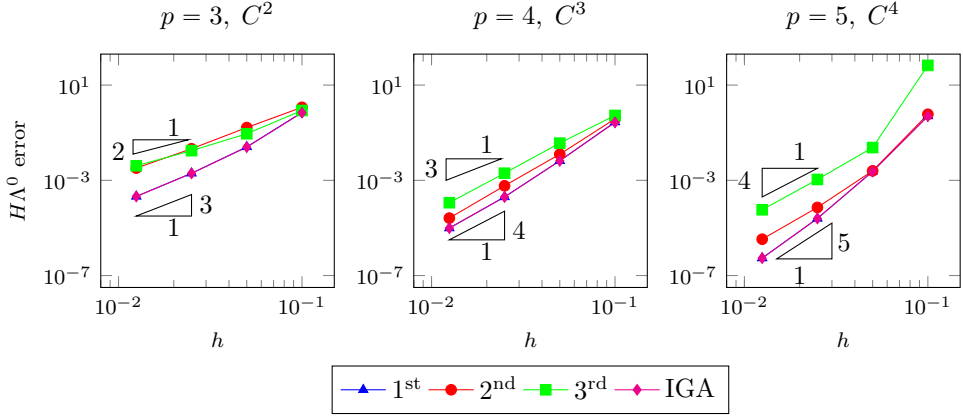
\begin{figure}[t!]
    \centering 
    \begin{tikzpicture}[]
\pgfplotsset{
    width=0.36\textwidth, 
    height=0.36\textwidth, 
    ymin=10^-7.5,ymax=200,xmin=8e-3,xmax=0.15, 
    ylabel near ticks,
    tick label style={font=\footnotesize},
    label style={font=\footnotesize},
    every x tick scale label/.style={at={(1,0)},anchor=north,yshift=-5pt,inner sep=0pt},
}
        \begin{groupplot}[ 
            group style={
            group size=3 by 1,
            horizontal sep=30pt
            },
        ]
            \nextgroupplot[ymode=log,xmode=log,title={$p=3,\;C^2$},ylabel={$H\Lambda^0$ error},xlabel={$h$},
            legend style={at={($(0,0)+(1cm,1cm)$)},legend columns=4,fill=none,draw=black,anchor=center,align=center},
            legend to name=fred]
            \addplot [blue, mark=triangle*, mark options={solid}] table [x index=0, y index=2, col sep=comma] {figures/laplace_results_mappedsquare_p3r2_global.csv};
            \addplot [red, mark=*, mark options={solid}] table [x index=0, y index=2, col sep=comma] {figures/laplace_results_mappedsquare_p3r2_l2.csv};
            \addplot [green, mark=square*, mark options={solid}] table [x index=0, y index=2, col sep=comma] {figures/laplace_results_mappedsquare_p3r2_llm.csv};
            \addplot [magenta, mark=diamond*, mark options={solid}] table [x index=0, y index=2, col sep=comma] {figures/laplace_results_mappedsquare_p3r2_standard_iga.csv};
            \draw (axis cs:0.024,2^2*10^-1.9)--(axis cs:0.012,10^-1.9)--(axis cs:0.012,2^2*10^-1.9) node[pos=0.05,anchor=north] {}--cyclenode[pos=0.0,anchor=west] {};
            \draw (axis cs:0.025,10^-4.5)--(axis cs:0.0125,10^-4.5)--(axis cs:0.025,2^3*10^-4.5) node[pos=0.05,anchor=north] {}--cyclenode[pos=0.0,anchor=west] {};
            \node[] at (axis cs:0.01,10^-1.8)   (c) {$2$};
            \node[] at (axis cs:0.03, 10^-4.0)   (c) {$3$};
            \node[] at (axis cs:0.02,10^-0.8)   (c) {$1$};
            \node[] at (axis cs:0.02,10^-5.0)   (c) {$1$};
            \addlegendentry{1$^{\text{st}}$};    
            \addlegendentry{2$^{\text{nd}}$};    
            \addlegendentry{3$^{\text{rd}}$};   
            \addlegendentry{IGA}; 
            \coordinate (c1) at (rel axis cs:0,1);
            \nextgroupplot[ymode=log,xmode=log,title={$p=4,\;C^3$},xlabel={$h$}]
            \addplot [blue, mark=triangle*, mark options={solid}] table [x index=0, y index=2, col sep=comma] {figures/laplace_results_mappedsquare_p4r3_global.csv};
            \addplot [red, mark=*, mark options={solid}] table [x index=0, y index=2, col sep=comma] {figures/laplace_results_mappedsquare_p4r3_l2.csv};
            \addplot [green, mark=square*, mark options={solid}] table [x index=0, y index=2, col sep=comma] {figures/laplace_results_mappedsquare_p4r3_llm.csv};
            \addplot [magenta, mark=diamond*, mark options={solid}] table [x index=0, y index=2, col sep=comma] {figures/laplace_results_mappedsquare_p4r3_standard_iga.csv};
            \draw (axis cs:0.024,2^3*10^-3)--(axis cs:0.012,10^-3)--(axis cs:0.012,2^3*10^-3) node[pos=0.05,anchor=north] {}--cyclenode[pos=0.0,anchor=west] {};
            \node[] at (axis cs:0.01,10^-2.5)   (c) {$3$};
            \draw (axis cs:0.025,10^-5.5)--(axis cs:0.0125,10^-5.5)--(axis cs:0.025,2^4*10^-5.5) node[pos=-0.05,anchor=north] {}--cyclenode[pos=1.0,anchor=west] {};
            \node[] at (axis cs:0.03, 0.9*10^-5.5+0.1*2^4*10^-5.5)   (b) {$4$};
            \node[] at (axis cs:0.02,10^-1.5)   (c) {$1$};
            \node[] at (axis cs:0.02,10^-6.0)   (c) {$1$};
            \nextgroupplot[ymode=log,xmode=log,title={$p=5,\;C^4$},xlabel={$h$}]
            \addplot [blue, mark=triangle*, mark options={solid}] table [x index=0, y index=2, col sep=comma] {figures/laplace_results_mappedsquare_p5r4_global.csv};
            \addplot [red, mark=*, mark options={solid}] table [x index=0, y index=2, col sep=comma] {figures/laplace_results_mappedsquare_p5r4_l2.csv};
            \addplot [green, mark=square*, mark options={solid}] table [x index=0, y index=2, col sep=comma] {figures/laplace_results_mappedsquare_p5r4_llm.csv};
            \addplot [magenta, mark=diamond*, mark options={solid}] table [x index=0, y index=2, col sep=comma] {figures/laplace_results_mappedsquare_p5r4_standard_iga.csv};
            \draw (axis cs:0.0125,2^4*10^-3.7)--(axis cs:0.0125,10^-3.7)--(axis cs:0.025,2^4*10^-3.7) node[pos=-0.05,anchor=north] {}--cyclenode[pos=1.0,anchor=west] {};
            \draw (axis cs:0.03,10^-6.3)--(axis cs:0.015,10^-6.3)--(axis cs:0.03,2^5 * 10^-6.3) node[pos=-0.05,anchor=north] {}--cyclenode[pos=1.0,anchor=west] {};
            \node[] at (axis cs:0.01, 10^-3.0)   (b) {$4$};
            \node[] at (axis cs:0.035, 0.95*10^-6.0+0.05*2^5*10^-6.0)   (c) {$5$};
            \node[] at (axis cs:0.02,10^-1.8)   (c) {$1$};
            \node[] at (axis cs:0.02,10^-7.0)   (c) {$1$};
            \coordinate (c2) at (rel axis cs:1,1);
        \end{groupplot}
        \coordinate (c3) at ($(c1)!.5!(c2)$);
        \node[below] at (c3 |- current bounding box.south)
        {\pgfplotslegendfromname{fred}};
    \end{tikzpicture} 
    \caption{Elliptic source problem for the parametric domain, with the nontrivial map of (\ref{eq:harmonic_map}): a comparison of the error in the $H\Lambda^0(\Omega)$ norm (or $H^1$ norm) for the three choices of discrete Hodge operators.}
    \label{fig:h1_error_noK}
\end{figure}

As a first test we show how the various approaches to construct a discrete Hodge operator work for a problem in which only one Hodge--star operator, namely $\star^1_\gamma$, is needed. The discretization of the static problem \eqref{eq:continuousproblem} \textcolor{black}{with $k=0$ fits the requirements, since taking $\gamma=1$ and $\beta=0$ gives us the standard Poisson problem}.
As a convergence test we choose the domain to be the usual unit square $\Omega = \hat{\Omega} := (0,1)^n$, but in the more interesting case in which the interior of the physical domain has been displaced by means of the map
\begin{equation}
x_i = u_i + \frac{a}{2}\, \textcolor{black}{\Pi_{j=1}^{n} \sin (2\pi(u_j-1/2))} , \; \text{ for } i \in \{ 1, 2, \dots , n\},\label{eq:harmonic_map}
\end{equation}
\noindent inspired by \cite{jain_construction_2020}, where $u_i$ are the parametric coordinates while $x_i$ are the physical ones on the mapped domain. The parameter $a$, with $0<a<1$, regulates the peak deformation in the mesh with respect to a Cartesian--orthogonal one, and we set $a = 0.3$. The effect of the above is to obtain a nontrivial map of $\hat{\Omega}$ to itself, shown in Fig.~\ref{fig:reparam}. In this way, we can evaluate how faithfully a non-polynomial perturbation of the metric is approximated by the Hodge--star operator on a test for which a well-known closed form solution is still easily computable.
For starting polynomial degrees $p=3,4,5$ in the primal sequence and various continuities we show performance of the standard Galerkin IGA, as implemented in \cite{Vazquez_2016aa}, and discrete Hodge--star $\star_h^1$ obtained either by \textcolor{black}{the global operator \eqref{eq:def_hodge_discrete1} (labeled $1^{st}$ in the plots), by the global $L^2$--projection onto the target space \eqref{eq:def_hodge_discrete2} (labeled $2^{nd}$ in the plots)}, or by the quasi-interpolant based on local projections discussed in Section~\ref{sec:dsplinecomplex} (labeled $3^{rd}$).
We show convergence rates in the ${H\Lambda^0(\Omega)}$ norm (or equivalently, the usual $H^1$ norm) in Fig.~\ref{fig:h1_error_noK} 
where asymptotic rates confirm the theoretical results of Section~\ref{sec:erroranalysis}.
Specifically, the first global operator converges as fast as the standard IGA Galerkin approach, as expected from the result in Theorem~\ref{thm:first-hodge}. The standard Galerkin IGA approach can be shown in fact to be algebraically derivable from the first approach, even though they are not computationally equivalent.
The convergence rate of the second operator 
and of the local operator usually have one order of convergence less and confirm the analysis in Theorem~\ref{thm:error_bound}.


Another interesting case in which a high order geometric method is desired is the case in which material parameters involved in the constitutive equations are not piecewise constants but instead smooth functions, This setting, in which we instead use the parametric domain $\Omega = \hat \Omega$, is particularly relevant in the wave propagation simulations, for instance in designing non reflective boundary conditions via perfectly matched layers~\cite{collino_optimizing_1998}.
We choose as an example the material parameter 
\begin{equation}\label{eq:smoothgamma}
\gamma = 1 + e^{\frac{(x-1/2)^2 + (y-1/2)^2}{20^2}},
\end{equation} 
and perform a new convergence test by means of a manufactured solution on the unit square with homogeneous boundary conditions. This second test differs from the first one since the material parameter acts as a weight in the projection but is not involved in the pullbacks of differential forms, differently from manipulating the geometry of the domain. \textcolor{black}{To compute the integrals we use the same quadrature rules as in the previous test, which are Gauss-Legendre rules with $(p+1)^2$ quadrature points per element. 
The convergence results in Fig.~\ref{fig:h1_error_smooth_coeff} confirm that we maintain the convergence rates of the previous case, in which we only manipulated the geometry, without adding more quadrature points.}

\begin{figure}[t!]
    \centering
    \begin{tikzpicture}[]
        \begin{loglogaxis}[width=0.45\columnwidth, height=0.45\columnwidth, ylabel={$\kappa$},xlabel={$h$},xmin=10^-2.0, ylabel near ticks,tick label style={font=\footnotesize},  legend style={font=\scriptsize, }, label style={font=\footnotesize},every x tick scale label/.style={at={(1,0)},anchor=north,yshift=-5pt,inner sep=0pt},legend pos= north east]
        \addplot [blue, mark=triangle*, mark options={solid}] table [x index=0, y index=3, col sep=comma] {figures/laplace_results_mappedsquare_p4r3_global.csv};
        \addplot [red, mark=*, mark options={solid}] table [x index=0, y index=3, col sep=comma] {figures/laplace_results_mappedsquare_p4r3_l2.csv};
        \addplot [green, mark=square*, mark options={solid}] table [x index=0, y index=3, col sep=comma] {figures/laplace_results_mappedsquare_p4r3_llm.csv};
        \addplot [magenta, mark=diamond*, mark options={solid}] table [x index=0, y index=3, col sep=comma] {figures/laplace_results_mappedsquare_p3r2_standard_iga.csv};
        \node[] at (axis cs:10^-1.6, 1)   (a) {$p=3,\;C^2$};
        \draw (axis cs:0.02,10^2.5)--(axis cs:0.04,10^2.5)--(axis cs:0.02,2^2.5*10^2.5) node[pos=0.05,anchor=north] {}--cyclenode[pos=0.0,anchor=west] {};
        \node[] at (axis cs:0.015, 10^3.0)   (c) {$2.5$};
        \node[] at (axis cs:0.03,10^2.2)   (c) {$1$};
        \legend{1$^{\text{st}}$ \\2$^{\text{nd}}$ \\3$^{\text{rd}}$ \\IGA\\}
        \end{loglogaxis}
    \end{tikzpicture}
    \begin{tikzpicture}[]
        \begin{loglogaxis}[width=0.45\columnwidth, height=0.45\columnwidth, ylabel={Number of nonzeros},xlabel={$h$}, xmin=10^-2.0, ylabel near ticks,tick label style={font=\footnotesize}, legend style={font=\scriptsize, }, label style={font=\footnotesize},every x tick scale label/.style={at={(1,0)},anchor=north,yshift=-5pt,inner sep=0pt},legend pos= north east]
        \addplot [blue, mark=triangle*, mark options={solid}] table [x index=0, y index=4, col sep=comma] {figures/laplace_results_mappedsquare_p3r2_global.csv};
        \addplot [red, mark=*, mark options={solid}] table [x index=0, y index=4, col sep=comma] {figures/laplace_results_mappedsquare_p3r2_l2.csv};
        \addplot [green, mark=square*, mark options={solid}] table [x index=0, y index=4, col sep=comma] {figures/laplace_results_mappedsquare_p3r2_llm.csv};
        \addplot [magenta, mark=diamond*, mark options={solid}] table [x index=0, y index=4, col sep=comma] {figures/laplace_results_mappedsquare_p3r2_standard_iga.csv};
        \node[] at (axis cs:10^-1.6, 10^-0.5)   (a) {$p=3,\;C^2$};
        \draw (axis cs:0.02,10^4.0)--(axis cs:0.04,10^4.0)--(axis cs:0.02,2^2.0*10^4.0) node[pos=0.05,anchor=north] {}--cyclenode[pos=0.0,anchor=west] {};
        \node[] at (axis cs:0.018, 10^4.2)   (c) {$2$};
        \node[] at (axis cs:0.03,10^3.8)   (c) {$1$};
        \legend{1$^{\text{st}}$ \\2$^{\text{nd}}$ \\3$^{\text{rd}}$ \\IGA\\}
        \end{loglogaxis}
    \end{tikzpicture} 
    \caption{Elliptic source problem for the parametric domain, with the nontrivial map of (\ref{eq:harmonic_map}): a comparison of condition number and sparsity degrade for the case of $h$ refinement for $p=3$ and $C^2$ splines in the first space of the primal sequence.}
    \label{fig:conditioning_and_sparsity_laplace}
\end{figure}
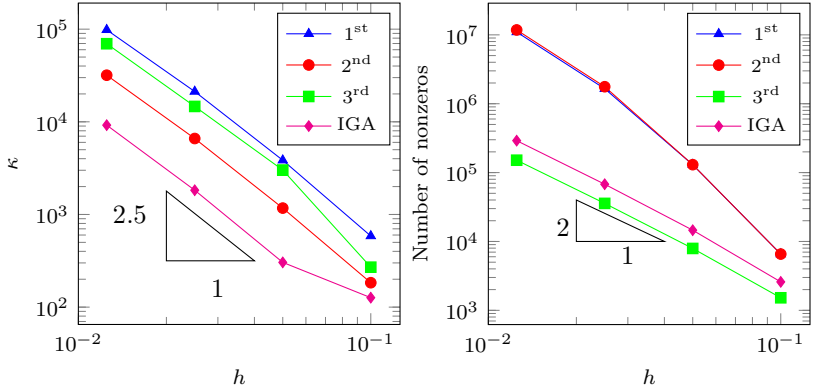

The condition number of the resulting system matrices for the test problem is studied in Fig.~\ref{fig:conditioning_and_sparsity_laplace}, left panel, \textcolor{black}{for the practically very relevant case of splines of degree $p=3$ and maximum smoothness for the $X^0$ space}.
When dealing with a square and invertible matrix $\mathbf{A}$, the condition number $\kappa(\mathbf{A})$ is defined in terms of the 2-norm of $\mathbf{A}$ as follows: 
\begin{equation*}
    \kappa(\mathbf{A}) = {\Vert\mathbf{A}\Vert_2} {\Vert\mathbf{A^{-1}}\Vert_2}.
\end{equation*}
In our examples $\kappa$ is computed for the matrix discretising the composition of exterior derivatives (pre-- and post--) with the $\star_h^1$ operator, \textcolor{black}{ that is, the matrix on the left-hand side of \eqref{eq:problem_matrix} with $\beta = 0$}. \textcolor{black}{The new approach in general yields the same asymptotic growth in condition number under mesh refinement as standard IGA, albeit degraded by a constant factor.} Furthermore, we remark that the system matrices \textcolor{black}{in \eqref{eq:problem_matrix}} are not in general symmetric and require an appropriate solving strategy (e.g. GMRES). \textcolor{black}{While this can be fixed for the two global operators by multiplying by the matrix $\tilde{\mathbf{K}}^k_{n-k}$ (which is equivalent to multiply by a test function, as we did for the analysis in Section~\ref{sec:erroranalysis}), the problem remains for the local operator based on quasi-interpolants. This is an important issue} when one aims to use the new machinery to evolve hyperbolic systems of equations in time, and it begs for further theoretical study.

The sparsity achieved in the final system matrix is also studied in Fig.~\ref{fig:conditioning_and_sparsity_laplace} on the right panel, \textcolor{black}{where we plot the number of nonzeros again for degree $p=3$ and maximum continuity of the spline spaces}. We see that the local projector yields sparser system matrices than the global operators, as expected, and even sparser than the standard IGA approach by a factor which is not dependent on $h$.

\begin{figure}[t!]
    \centering 
    \begin{tikzpicture}[]
        \pgfplotsset{
            width=0.36\textwidth, 
            height=0.36\textwidth, 
            ymin=10^-7.5,ymax=5*10^0,xmin=8e-3,xmax=0.15, 
            ylabel near ticks,
            tick label style={font=\footnotesize},
            label style={font=\footnotesize},
            every x tick scale label/.style={at={(1,0)},anchor=north,yshift=-5pt,inner sep=0pt},
        }
        \begin{groupplot}[ 
            group style={
            group size=3 by 1,
            horizontal sep=30pt
            },
        ]
            \nextgroupplot[ymode=log,xmode=log,title={$p=2,\;C^1$},ylabel={$H\Lambda^0$ error},xlabel={$h$},
            legend style={at={($(0,0)+(1cm,1cm)$)},legend columns=4,fill=none,draw=black,anchor=center,align=center},
            legend to name=barney]
            \addplot [blue, mark=triangle*, mark options={solid}] table [x index=0, y index=2, col sep=comma] {figures/laplace_results_square_highreg_smoothgamma_p2r1_global.csv};
            \addplot [red, mark=*, mark options={solid}] table [x index=0, y index=2, col sep=comma] {figures/laplace_results_square_highreg_smoothgamma_p2r1_l2.csv};
            \addplot [green, mark=square*, mark options={solid}] table [x index=0, y index=2, col sep=comma] {figures/laplace_results_square_highreg_smoothgamma_p2r1_llm.csv};
            \addplot [magenta, mark=diamond*, mark options={solid}] table [x index=0, y index=2, col sep=comma] {figures/laplace_results_square_highreg_smoothgamma_p2r1_standard_iga.csv};
            \draw (axis cs:0.012,0.5*10^-0.4)--(axis cs:0.024,10^-0.4)--(axis cs:0.012,10^-0.4) node[pos=0.05,anchor=north] {}--cyclenode[pos=0.0,anchor=west] {};
            \node[] at (axis cs:0.01, 10^-0.5)   (c) {$1$};
            \draw (axis cs:0.024,2^2*10^-3.4)--(axis cs:0.024,10^-3.4)--(axis cs:0.012,10^-3.4) node[pos=0.05,anchor=north] {}--cyclenode[pos=0.0,anchor=west] {};
            \node[] at (axis cs:0.03, 10^-3.2)   (c) {$2$};
            \node[] at (axis cs:0.02,10^-0.0125)   (c) {$1$};
            \node[] at (axis cs:0.02,10^-4.0)   (c) {$1$};
            \addlegendentry{1$^{\text{st}}$};    
            \addlegendentry{2$^{\text{nd}}$};    
            \addlegendentry{3$^{\text{rd}}$};   
            \addlegendentry{IGA}; 
            \coordinate (c1) at (rel axis cs:0,1);
            \nextgroupplot[ymode=log,xmode=log,title={$p=3,\;C^2$},xlabel={$h$}]
            \addplot [blue, mark=triangle*, mark options={solid}] table [x index=0, y index=2, col sep=comma] {figures/laplace_results_square_highreg_smoothgamma_p3r2_global.csv};
            \addplot [red, mark=*, mark options={solid}] table [x index=0, y index=2, col sep=comma] {figures/laplace_results_square_highreg_smoothgamma_p3r2_l2.csv};
            \addplot [green, mark=square*, mark options={solid}] table [x index=0, y index=2, col sep=comma] {figures/laplace_results_square_highreg_smoothgamma_p3r2_llm.csv};
            \addplot [magenta, mark=diamond*, mark options={solid}] table [x index=0, y index=2, col sep=comma] {figures/laplace_results_square_highreg_smoothgamma_p3r2_standard_iga.csv};
            \draw (axis cs:0.024,2^2*10^-2.4)--(axis cs:0.012,10^-2.4)--(axis cs:0.012,2^2*10^-2.4) node[pos=0.05,anchor=north] {}--cyclenode[pos=0.0,anchor=west] {};
            \draw (axis cs:0.025,10^-5.0)--(axis cs:0.0125,10^-5.0)--(axis cs:0.025,2^3*10^-5.0) node[pos=0.05,anchor=north] {}--cyclenode[pos=0.0,anchor=west] {};
            \node[] at (axis cs:0.01, 10^-2.2)   (c) {$2$};
            \node[] at (axis cs:0.03, 10^-4.5)   (c) {$3$};
            \node[] at (axis cs:0.02,10^-1.2)   (c) {$1$};
            \node[] at (axis cs:0.02,10^-5.5)   (c) {$1$};
            \nextgroupplot[ymode=log,xmode=log,title={$p=4,\;C^3$},xlabel={$h$}]
            \addplot [blue, mark=triangle*, mark options={solid}] table [x index=0, y index=2, col sep=comma] {figures/laplace_results_square_highreg_smoothgamma_p4r3_global.csv};
            \addplot [red, mark=*, mark options={solid}] table [x index=0, y index=2, col sep=comma] {figures/laplace_results_square_highreg_smoothgamma_p4r3_l2.csv};
            \addplot [green, mark=square*, mark options={solid}] table [x index=0, y index=2, col sep=comma] {figures/laplace_results_square_highreg_smoothgamma_p4r3_llm.csv};
            \addplot [magenta, mark=diamond*, mark options={solid}] table [x index=0, y index=2, col sep=comma] {figures/laplace_results_square_highreg_smoothgamma_p4r3_standard_iga.csv};
            \draw (axis cs:0.024,2^3*10^-3.4)--(axis cs:0.012,10^-3.4)--(axis cs:0.012,2^3*10^-3.4) node[pos=0.05,anchor=north] {}--cyclenode[pos=0.0,anchor=west] {};
            \node[] at (axis cs:0.01,10^-3.0)   (c) {$3$};
            \draw (axis cs:0.025,10^-6.4)--(axis cs:0.0125,10^-6.4)--(axis cs:0.025,2^4*10^-6.4) node[pos=-0.05,anchor=north] {}--cyclenode[pos=1.0,anchor=west] {};
            \node[] at (axis cs:0.03,10^-6.0)   (b) {$4$};
            \node[] at (axis cs:0.02,10^-2.0)   (c) {$1$};
            \node[] at (axis cs:0.02,10^-7.0)   (c) {$1$};
            \coordinate (c2) at (rel axis cs:1,1);
        \end{groupplot}
        \coordinate (c3) at ($(c1)!.5!(c2)$);
        \node[below] at (c3 |- current bounding box.south)
        {\pgfplotslegendfromname{barney}};
    \end{tikzpicture} 
    \caption{Static problem in the case of the parametric domain, with the smooth perturbation of the material parameter $\gamma$ introduced in \eqref{eq:smoothgamma}: a comparison of the error in the $H\Lambda^0(\Omega)$ norm (or $H^1$ norm) for various polynomial degrees for the three choices of discrete Hodge operators.}
    \label{fig:h1_error_smooth_coeff}
\end{figure}
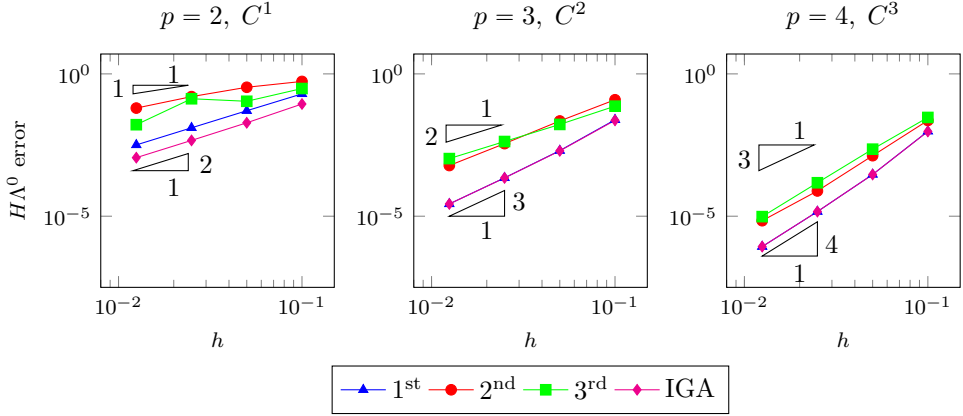

\subsection{Second order eigenvalue problem}

In the case of the eigenvalue problem \eqref{eq:evp_with_forms} \textcolor{black}{we are particularly interested in Maxwell eigenvalue problem, for which we set $k=1$ and the material parameters $\gamma = \beta = 1$. The} most important property we want to investigate in the proposed approach is its spectral accuracy. This test is fundamental in a twofold sense: firstly, the absence of numerical spurious eigenfunctions in the kernel of the exterior derivative operators is a key feature of structure preserving methods, secondly the optimal $2p$ convergence rate in eigenvalue errors is sought to confirm that the approach is indeed a high-order one.

\begin{figure}[t!]
    \centering 
    \begin{tikzpicture}[]
            \begin{axis}[width=0.45\columnwidth, height=0.45\columnwidth, ylabel={},xlabel={Eigenvalue index}, 
                ymin=0,ymax=110,xmin=0,xmax=21, ylabel={$\tau$}, ylabel near ticks,tick label style={font=\footnotesize},  legend style={font=\scriptsize, }, label style={font=\footnotesize},every x tick scale label/.style={at={(1,0)},anchor=north,yshift=-5pt,inner sep=0pt},legend pos= south east]
            \addplot [red, only marks, mark=diamond*, 
            ] table [x index=0, y index=2, col sep=comma] {figures/maxwell_results_eighthring_reference_p4r3_refinement1.csv};
            \addplot [green, only marks, mark=x, 
            ] table [x index=0, y index=3, col sep=comma] {figures/maxwell_results_eighthring_reference_p4r3_refinement1.csv};
            \addplot [black, only marks, mark=square, mark options={}
            ] table [x index=0, y index=4, col sep=comma] {figures/maxwell_results_eighthring_reference_p4r3_refinement1.csv};
            \node[] at (axis cs:10,99)   (b) {$1/8$ of a ring};
            \legend{2$^{\text{nd}}$ \\3$^{\text{rd}}$ \\Reference\\}
            \end{axis}
    \end{tikzpicture}
    \begin{tikzpicture}[]
            \begin{axis}[width=0.45\columnwidth, height=0.45\columnwidth, ylabel={$\tau/\pi^2$},xlabel={Eigenvalue index}, 
                ymin=0,ymax=9,xmin=0,xmax=21, ylabel near ticks,tick label style={font=\footnotesize},  legend style={font=\scriptsize, }, label style={font=\footnotesize},every x tick scale label/.style={at={(1,0)},anchor=north,yshift=-5pt,inner sep=0pt},legend pos= south east]
            \addplot [red, only marks, mark=diamond*, 
            ] table [x index=0, y index=3, col sep=comma] {figures/maxwell_results_cube_p4r3_refinement3.csv};
            \addplot [green, only marks, mark=x, 
            ] table [x index=0, y index=4, col sep=comma] {figures/maxwell_results_cube_p4r3_refinement3.csv};
            \addplot [black, only marks, mark=square, mark options={}
            ] table [x index=0, y index=1, col sep=comma] {figures/maxwell_results_cube_p4r3_refinement3.csv};
            \node[] at (axis cs:10,8)   (b) {Unit cube};
            \legend{2$^{\text{nd}}$ \\3$^{\text{rd}}$ \\Reference\\}
            \end{axis}
    \end{tikzpicture}
    \caption{Spectral fidelity: \textcolor{black}{comparison of the computed eigenvalues for the second and third discrete Hodge operators with reference eigenvalues for one eighth of a ring, and for a deformation of the unit cube.}}
    \label{fig:ev_correctness}
\end{figure}
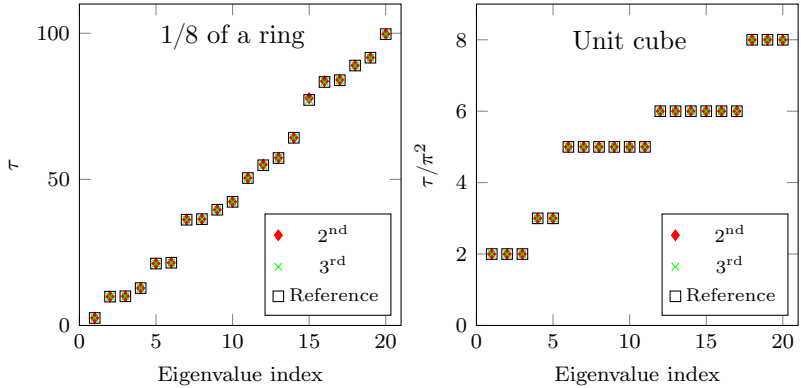

\begin{figure}[t!]
    \centering 
    \begin{tikzpicture}[]
            \begin{loglogaxis}[width=0.45\columnwidth, title={$p=4,\;C^3$}, height=0.45\columnwidth, ylabel={Eigenvalue error},xlabel={$h$}, 
                ymin=10^-12.2,ymax=10^-0.2,xmin=0.02,xmax=2, ylabel near ticks,tick label style={font=\footnotesize},  legend style={font=\scriptsize, }, label style={font=\footnotesize},every x tick scale label/.style={at={(1,0)},anchor=north,yshift=-5pt,inner sep=0pt},legend pos= south east]
            \addplot [blue, mark=triangle*, mark options={solid}] table [x index=0, y index=20, col sep=comma] {figures/maxwell_results_cube_p4r3_global.csv};
            \addplot [red, mark=*, mark options={solid}] table [x index=0, y index=20, col sep=comma] {figures/maxwell_results_cube_p4r3_l2.csv};
            \addplot [green, mark=square*, mark options={solid}] table [x index=0, y index=20, col sep=comma] {figures/maxwell_results_cube_p4r3_llm.csv};
            \addplot [magenta, mark=diamond*, mark options={solid}] table [x index=0, y index=20, col sep=comma] {figures/maxwell_results_cube_p4r3_standard_iga.csv};
            \draw (axis cs:0.05,2^4*10^-3.7)--(axis cs:0.05,10^-3.7)--(axis cs:0.1,2^4*10^-3.7) node[pos=-0.05,anchor=north] {}--cyclenode[pos=1.0,anchor=west] {};
            \draw (axis cs:0.05,2^6 * 10^-7.0)--(axis cs:0.05,10^-7.0)--(axis cs:0.1,2^6 * 10^-7.0) node[pos=-0.05,anchor=north] {}--cyclenode[pos=1.0,anchor=west] {};
            \draw (axis cs:0.2,10^-8.7)--(axis cs:0.1,10^-8.7)--(axis cs:0.2,2^8 * 10^-8.7) node[pos=-0.05,anchor=north] {}--cyclenode[pos=1.0,anchor=west] {};
            \node[] at (axis cs:0.04, 0.95*10^-3.0+0.05*2^4*10^-3.0)   (b) {$4$};
            \node[] at (axis cs:0.04, 0.9*10^-7.0+0.1*2^6*10^-7.0)   (c) {$6$};
            \node[] at (axis cs:0.25, 0.9*10^-9.0+0.1*2^8*10^-9.0)   (c) {$8$};
            \node[] at (axis cs:0.07, 10^-2.0)   (b) {$1$};
            \node[] at (axis cs:0.07, 10^-4.5)   (c) {$1$};
            \node[] at (axis cs:0.15, 10^-9.5)   (c) {$1$};
            \legend{1$^{\text{st}}$ \\2$^{\text{nd}}$ \\3$^{\text{rd}}$ \\IGA\\}
            \end{loglogaxis}
        \end{tikzpicture}
        \begin{tikzpicture}[]
                \begin{loglogaxis}[width=0.45\columnwidth, title={$p=5,\;C^4$}, height=0.45\columnwidth, ylabel={},xlabel={$h$}, 
                    ymin=10^-12.2,ymax=10^-0.2,xmin=0.02,xmax=2, ylabel near ticks,tick label style={font=\footnotesize},  legend style={font=\scriptsize, }, label style={font=\footnotesize},every x tick scale label/.style={at={(1,0)},anchor=north,yshift=-5pt,inner sep=0pt},legend pos= south east]
                \addplot [blue, mark=triangle*, mark options={solid}] table [x index=0, y index=20, col sep=comma] {figures/maxwell_results_cube_p5r4_global.csv};
                \addplot [red, mark=*, mark options={solid}] table [x index=0, y index=20, col sep=comma] {figures/maxwell_results_cube_p5r4_l2.csv};
                \addplot [green, mark=square*, mark options={solid}] table [x index=0, y index=20, col sep=comma] {figures/maxwell_results_cube_p5r4_llm.csv};
                \addplot [magenta, mark=diamond*, mark options={solid}] table [x index=0, y index=20, col sep=comma] {figures/maxwell_results_mappedcube_p5r4_standard_iga.csv};
                \draw (axis cs:0.05,2^4*10^-3.7)--(axis cs:0.05,10^-3.7)--(axis cs:0.1,2^4*10^-3.7) node[pos=-0.05,anchor=north] {}--cyclenode[pos=1.0,anchor=west] {};
                \draw (axis cs:0.08,2^8 * 10^-7.5)--(axis cs:0.08,10^-7.5)--(axis cs:0.16,2^8 * 10^-7.5) node[pos=-0.05,anchor=north] {}--cyclenode[pos=1.0,anchor=west] {};
                \draw (axis cs:0.16,10^-11.0)--(axis cs:0.08,10^-11.0)--(axis cs:0.16,2^10 * 10^-11.0) node[pos=-0.05,anchor=north] {}--cyclenode[pos=1.0,anchor=west] {};
                \node[] at (axis cs:0.04, 0.95*10^-3.2+0.05*2^4*10^-3.2)   (b) {$4$};
                \node[] at (axis cs:0.06, 0.9*10^-7.0+0.1*2^6*10^-7.0)   (c) {$8$};
                \node[] at (axis cs:0.2, 10^-10.0)   (c) {$10$};
                \node[] at (axis cs:0.07, 10^-2.0)   (b) {$1$};
                \node[] at (axis cs:0.12, 10^-4.5)   (c) {$1$};
                \node[] at (axis cs:0.12, 10^-11.5)   (c) {$1$};
                \legend{1$^{\text{st}}$ \\2$^{\text{nd}}$ \\3$^{\text{rd}}$ \\IGA\\}
                \end{loglogaxis}
            \end{tikzpicture}
    \caption{Eigenvalue problem: convergence in the twentieth eigenvalue error for the various choices of projectors for the unit cube in three dimensions. The case of starting degrees $p=4$, $p=5$ with maximum continuity of spline spaces.}
    \label{fig:ev_error}
\end{figure}
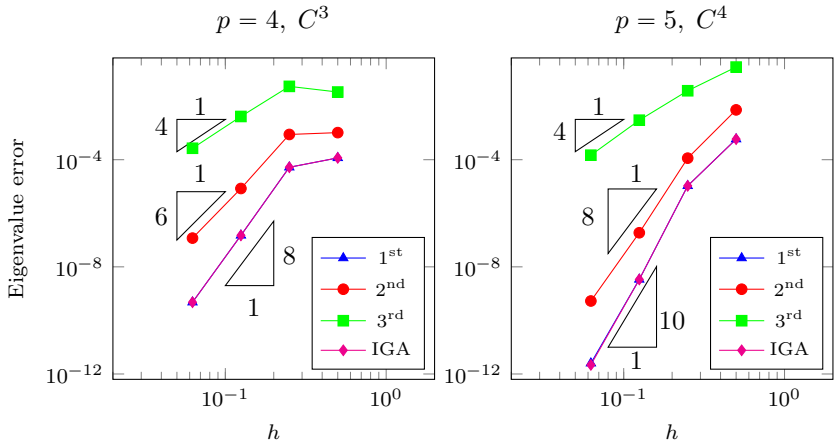

The first property is easily proved for the first global discrete Hodge--star operator proposed through its algebraic equivalence with a standard Galerkin formulation. Furthermore it is also reasonably expected to hold for the second kind of Hodge--star operator, being it equivalent to an $L^2$ projection, combining insight from~\cite{BRSV11} on the existence of projectors commuting with the De-Rham sequence operators and results in~\cite{Buffa_2020aa} about the stable $L^2$ pairing of univariate spline spaces of degree $(p,p-2)$. The fact that spectral correctness is observed in the case of the third proposed operator is instead considerably less trivial and deserves further attention, since the proposed local projector does not commute with differential operators in the sequence.
A full proof of this statement goes beyond the scope of the present article. 
Nevertheless, we show that this property is numerically corroborated in Fig.~\ref{fig:ev_correctness} both in two dimensions for the case of the domain one eighth of a ring and in three dimensions for the usual unit cube, albeit with the nontrivial map already used for the square in the problem of Section~\ref{sec:elliptic_numerical}, with the same value of $a=0.3$. \textcolor{black}{In the first case the reference eigenvalues are computed by a Galerkin method in a very fine mesh, while in the second case the exact eigenvalues are known.} 
The convergence rate in eigenvalue errors, computed for the second eigenvalue in the case of the eighth of a ring cavity and for the twentieth eigenvalue in the case of the cube and shown in Fig.~\ref{fig:ev_error_8th_ring} and Fig.~\ref{fig:ev_error}, respectively, are again the optimal ones for the approach which is algebraically equivalent to the standard Galerkin IGA one. 
The other two suggested projectors instead lose orders of convergence upon $h-$refinement. Since we perform $L^2$ projections into spaces of splines with mixed degree between $p-1$ and $p-2$ the $2^\text{nd}$ operator actually exhibits the \textcolor{black}{expected convergence rate deriving from Theorem~\ref{thm:error_bound}}, i.e. $h^{2(p-1)}$. The local operator instead loses more (and the behaviour gets more evident for increasing polynomial degrees) as a result of the combination of loss of symmetry of the system matrix and the filtering of tails of some basis functions in the local projection operator. The high order nature of the method is still preserved, yet future work is needed for optimising the presented approach for eigenvalue solvers.



\begin{figure}[t!]
    \centering 
    \begin{tikzpicture}[]
    \begin{loglogaxis}[width=0.45\columnwidth, height=0.45\columnwidth, ylabel={Eigenvalue error},xlabel={$h$}, title={$p=3,\;C^2$},
    ymin=10^-10,ymax=10^-1,xmin=0.01,xmax=0.5, ylabel near ticks,tick label style={font=\footnotesize},  legend style={font=\scriptsize, }, label style={font=\footnotesize},every x tick scale label/.style={at={(1,0)},anchor=north,yshift=-5pt,inner sep=0pt},legend pos= south east]
    \addplot [blue, mark=triangle*, mark options={solid}] table [x index=0, y index=2, col sep=comma] {figures/maxwell_results_eighthring_href_p3r2_global.csv};
    \addplot [red, mark=*, mark options={solid}] table [x index=0, y index=2, col sep=comma] {figures/maxwell_results_eighthring_href_p3r2_l2.csv};
    \addplot [green, mark=square*, mark options={solid}] table [x index=0, y index=2, col sep=comma] {figures/maxwell_results_eighthring_href_p3r2_llm.csv};
    \addplot [magenta, mark=diamond*, mark options={solid}] table [x index=0, y index=2, col sep=comma] {figures/maxwell_results_eighthring_href_p3r2_standard_iga.csv};
    \draw (axis cs:0.05,2^4*10^-3.2)--(axis cs:0.05,10^-3.2)--(axis cs:0.1,2^4*10^-3.2) node[pos=-0.05,anchor=north] {}--cyclenode[pos=1.0,anchor=west] {};
    \draw (axis cs:0.02,2^6 * 10^-9.0)--(axis cs:0.02,10^-9.0)--(axis cs:0.04,2^6 * 10^-9.0) node[pos=-0.05,anchor=north] {}--cyclenode[pos=1.0,anchor=west] {};
    \node[] at (axis cs:0.04, 10^-2.5)   (b) {$4$};
    \node[] at (axis cs:0.015, 0.9*10^-9.0+0.1*2^6*10^-9.0)   (c) {$6$};
    \node[] at (axis cs:0.07, 10^-1.5)   (b) {$1$};
    \node[] at (axis cs:0.03, 10^-6.8)   (c) {$1$};
    \legend{1$^{\text{st}}$ \\2$^{\text{nd}}$ \\3$^{\text{rd}}$ \\IGA\\}
    \end{loglogaxis}
    \end{tikzpicture}
    \begin{tikzpicture}[]
    \begin{loglogaxis}[width=0.45\columnwidth, height=0.45\columnwidth, ylabel={},xlabel={$h$}, title={$p=3,\;C^1$},
    ymin=10^-10,ymax=10^-1,xmin=0.01,xmax=0.5, ylabel near ticks,tick label style={font=\footnotesize},  legend style={font=\scriptsize, }, label style={font=\footnotesize},every x tick scale label/.style={at={(1,0)},anchor=north,yshift=-5pt,inner sep=0pt},legend pos= south east]
    \addplot [blue, mark=triangle*, mark options={solid}] table [x index=0, y index=2, col sep=comma] {figures/maxwell_results_eighthring_href_p3r1_global.csv};
    \addplot [red, mark=*, mark options={solid}] table [x index=0, y index=2, col sep=comma] {figures/maxwell_results_eighthring_href_p3r1_l2.csv};
    \addplot [green, mark=square*, mark options={solid}] table [x index=0, y index=2, col sep=comma] {figures/maxwell_results_eighthring_href_p3r1_llm.csv};
    \addplot [magenta, mark=diamond*, mark options={solid}] table [x index=0, y index=2, col sep=comma] {figures/maxwell_results_eighthring_href_p3r1_standard_iga.csv};
    \draw (axis cs:0.06,2^4*10^-3.7)--(axis cs:0.06,10^-3.7)--(axis cs:0.12,2^4*10^-3.7) node[pos=-0.05,anchor=north] {}--cyclenode[pos=1.0,anchor=west] {};
    \draw (axis cs:0.02,2^6 * 10^-9.0)--(axis cs:0.02,10^-9.0)--(axis cs:0.04,2^6 * 10^-9.0) node[pos=-0.05,anchor=north] {}--cyclenode[pos=1.0,anchor=west] {};
    \node[] at (axis cs:0.05, 0.95*10^-3.2+0.05*2^4*10^-3.2)   (b) {$4$};
    \node[] at (axis cs:0.015, 0.9*10^-9.0+0.1*2^6*10^-9.0)   (c) {$6$};
    \node[] at (axis cs:0.08, 10^-2.0)   (b) {$1$};
    \node[] at (axis cs:0.03, 10^-6.8)   (c) {$1$};
    \legend{1$^{\text{st}}$ \\2$^{\text{nd}}$ \\3$^{\text{rd}}$ \\IGA\\}
    \end{loglogaxis}
    \end{tikzpicture}
    \caption{Eigenvalue problem: convergence in the second eigenvalue error for the various choices of projectors for the eighth ring in two dimensions. The case of starting degree $p=3$ both with $C^2$ and $C^1$ continuous spline spaces.}
    \label{fig:ev_error_8th_ring}
\end{figure}
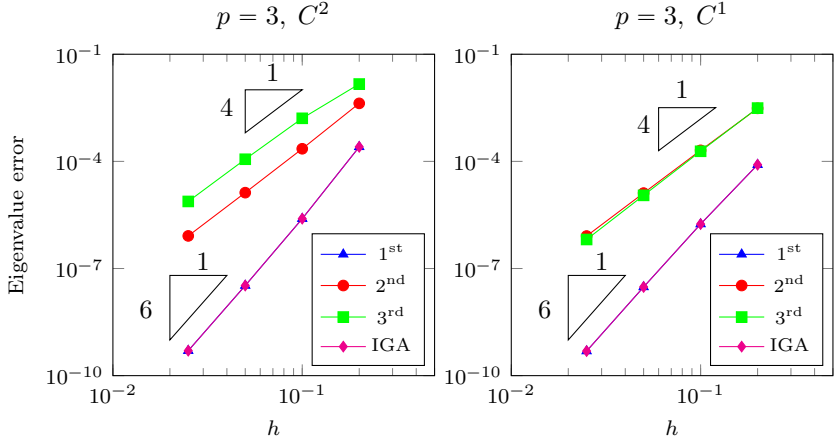

\section{Conclusions and future work} \label{Sec:conclusion}

We have introduced a new kind of geometric method based on spline spaces satisfying a de Rham complex. Thanks to the high continuity of splines, the dual complex is defined with a simple change of degree, and without explicitly constructing a dual mesh. Since B-splines are defined by tensor-product, the method works in arbitrary dimension, and the exterior derivative is always given by incidence matrices of a Cartesian mesh, also for curved domains. The method must be completed with discrete Hodge--star operators that link the primal and the dual complex, and we have introduced three different operators: the first two are global, while the third one is local and gives a reduced sparsity pattern. We have analyzed the discretization error and the order of the method applied to an elliptic problem, for the three different operators, following a variational approach. We proved that using the first operator our method is equivalent to the standard Galerkin method used in IGA, for which the same order of convergence is attained, while for the other two operators the expected order of convergence is reduced by one, due to the reduced degree of the splines in the dual complex. 

In a forthcoming paper we will analyze the method for the solution of the time domain Maxwell's equations. So far, the method is restricted to the case of a single patch, i.e., the domain is the image of the parametric domain through a given parameterization. Our goal is to extend it to geometries formed by multiple patches, the difficulty being that the regularity across patches is the same as for standard finite elements, and therefore the dual complex cannot be simply constructed by reducing the degree. The same issue arises when the domain is formed by multiple materials, and we expect that the necessary modifications will be the same for both cases.


\begin{appendices}
\end{appendices}

\bibliographystyle{siam}
\bibliography{biblio}

\end{document}